\documentclass[11pt]{amsart}

\usepackage{amsxtra}
\usepackage{amsfonts}
\usepackage[latin1]{inputenc}
\usepackage{epsfig}
\usepackage[T1]{fontenc}
\usepackage{amsmath,amssymb,latexsym,amsthm}
\usepackage[all]{xy}
\usepackage{mathrsfs}
\usepackage{psfrag}

\newtheorem{theorem}{Theorem}[section]
\newtheorem{lemma}[theorem]{Lemma}
\newtheorem{corollary}[theorem]{Corollary}
\newtheorem{definition}[theorem]{Definition}
\newtheorem{proposition}[theorem]{Proposition}
\newtheorem{conjecture}[theorem]{Conjecture}
\newtheorem{remark}[theorem]{Remark}
\newtheorem{example}[theorem]{Example}

\newtheorem{observation}[theorem]{Observation}

\numberwithin{equation}{section}

\CompileMatrices

\begin{document}

\title{Topological Hochschild homology and cohomology of $A_\infty$ ring spectra}

\author{Vigleik Angeltveit} \thanks{This research was partially conducted during the period the author was employed by the Clay Mathematics Institute as a Liftoff Fellow}

\maketitle

\newcommand{\A}{\mathcal{A}}
\newcommand{\B}{\mathcal{B}}
\newcommand{\C}{\mathcal{C}}
\newcommand{\D}{\mathcal{D}}
\newcommand{\E}{\mathcal{E}}
\newcommand{\F}{\mathbb{F}}
\newcommand{\I}{\mathcal{I}}
\newcommand{\K}{\mathcal{K}}
\newcommand{\M}{\mathcal{M}}
\newcommand{\R}{\mathbb{R}}
\newcommand{\V}{\mathcal{V}}
\newcommand{\W}{\mathbb{W}}
\newcommand{\Z}{\mathbb{Z}}
\newcommand{\sma}{\wedge}
\newcommand{\lar}{\longrightarrow}
\newcommand{\sar}{\rightarrow}
\newcommand{\colim}{\varinjlim}
\newcommand{\hEn}{\widehat{E(n)}}
\newcommand{\bD}{\mathbf{\Delta}}
\newcommand{\noco}{\mathbf{\Delta} \Sigma}
\newcommand{\Top}{\mathcal{T}\!\mathit{op}}

\def\arrow#1{\overset{#1}{\lar}}

\begin{abstract}
Let $A$ be an $A_\infty$ ring spectrum. We use the description from \cite{An_cyclic} of the cyclic bar and cobar construction to give a direct definition of topological Hochschild homology and cohomology of $A$ using the Stasheff associahedra and another family of polyhedra called cyclohedra. This construction builds the maps making up the $A_\infty$ structure into $THH(A)$, and allows us to study how $THH(A)$ varies over the moduli space of $A_\infty$ structures on $A$.

As an example, we study how topological Hochschild cohomology of Morava $K$-theory varies over the moduli space of $A_\infty$ structures and show that in the generic case, when a certain matrix describing the noncommutativity of the multiplication is invertible, topological Hochschild cohomology of $2$-periodic Morava $K$-theory is the corresponding Morava $E$-theory. If the $A_\infty$ structure is ``more commutative'', topological Hochschild cohomology of Morava $K$-theory is some extension of Morava $E$-theory.
\end{abstract}

\section{Introduction}
The main goal of this paper is to calculate topological Hochschild homology and cohomology of $A_\infty$ ring spectra such as Morava $K$-theory, and because $THH$ is sensitive to the $A_\infty$ structure we need to study the set (or space) of $A_\infty$ structures on a spectrum more closely. In particular, $THH(A)$ is sensitive to whether or not the multiplication is commutative, which is not so surprising if we think of topological Hochschild cohomology of $A$ as a version of the center of $A$. The Morava $K$-theory spectra are not even homotopy commutative at $p=2$, and at odd primes there is something noncommutative about the $A_p$ structure. Moreover, if we make Morava $K$-theory $2$-periodic it has many different homotopy classes of homotopy associative multiplications, most of which are noncommutative and all of which can be extended to $A_\infty$ structures.

Let us write $THH(A)$ for either topological Hochschild homology or cohomology of $A$, while using $THH^S(A)$ for topological Hochschild homology and $THH_S(A)$ for topological Hochschild cohomology. While $THH$ of an $A_\infty$ ring spectrum $A$ can be defined in a standard way after replacing it with a weakly equivalent $S$-algebra $\tilde{A}$, it is hard to see how the $A_\infty$ structure on $A$ affects the $S$-algebra structure on $\tilde{A}$. Instead we will define $THH(A)$ directly in terms of associahedra and cyclohedra, following \cite{An_cyclic}. In this way, the maps making up the $A_\infty$ structure on $A$ play a direct role, instead of being hidden away in the construction of a strictly associative replacement of $A$.

Our construction has the added advantage that given an $A_n$ structure on $A$, we can define spectra $sk_{n-1} THH^S(A)$ and $Tot^{n-1} THH_S(A)$. If the $A_n$ structure can be extended to an $A_\infty$ structure, $sk_{n-1} THH^S(A)$ coincides with the $(n-1)$-skeleton of $THH^S(A)$, but it is defined even if the $A_n$ structure cannot be extended, and similarly for $Tot^{n-1} THH_S(A)$.

We will be especially interested in $THH_R(R/I)$ and $THH^R(R/I)$, where $R$ is an even commutative $S$-algebra and $I=(x_1,\ldots,x_m)$ is a regular ideal. In this case any homotopy associative multiplication on $R/I$ can be extended to an $A_\infty$ structure (Corollary \ref{cor:evenquotientainfty}) and we have spectral sequences
\begin{eqnarray}
E_2^{*,*}=R_*/I[q_1,\ldots,q_m] \Longrightarrow  \pi_* THH_R(R/I); \label{eq:THHR/ISSintro1} \\
E^2_{*,*}=\Gamma_{R_*/I}[\bar{q}_1,\ldots,\bar{q}_m] \Longrightarrow  \pi_* THH^R(R/I). \label{eq:THHR/ISSintro2}
\end{eqnarray}
Here $\Gamma_{R_*/I}[\bar{q}_1,\ldots,\bar{q}_m]$ denotes a divided power algebra, though topological Hochschild homology is not generally a ring spectrum, so this has to be interpreted additively only. Topological Hochschild cohomology on the other hand is a ring spectrum. By the Deligne conjecture (see for example \cite{McSm02}) topological Hochschild cohomology admits an action of an $E_2$ operad, and in particular $THH_R(A)$ is a homotopy commutative ring spectrum. The first spectral sequence is a spectral sequence of commutative algebras, and it acts on the second spectral sequence in a natural way with $q_i \gamma_j(\bar{q}_i)=\gamma_{j-1}(\bar{q}_i)$, corresponding to the natural action $THH_R(A) \sma_R THH^R(A) \sar THH^R(A)$.

The above spectral sequences collapse, because they are concentrated in even total degree. But there are hidden extensions, and we can find these extensions by studying the $A_\infty$ structure on $R/I$ more closely. In the easiest case, when $I=(x)$, there is a hidden extension of height $n-1$ if the $A_n$ structure is ``noncommutative'' in a sense we will make precise in section \ref{sec:cyclic}. In general a careful analysis of the $A_\infty$ structure gives all of the extensions, and shows that
\begin{equation}
\pi_* THH_R(R/I) \cong R_*[[q_1,\ldots,q_m]]/(x_1-f_1,\ldots,x_m-f_m)
\end{equation}
for power series $f_1,\ldots,f_m$ which depend on the $A_\infty$ structure on $R/I$. By varying the $A_\infty$ structure the power series $f_i$ change, and this shows how $THH_R(R/I)$ varies over the moduli space of $A_\infty$ structures on $R/I$.

\subsection*{Organization}
In section \ref{sec:THHdef} we define topological Hochschild homology and cohomology using associahedra and cyclohedra following \cite{An_cyclic}, and prove that our definition agrees with the standard definition whenever they overlap. In section \ref{sec:Ainftyobtheory} we improve Robinson's obstruction theory \cite{Ro88} for endowing a spectrum with an $A_\infty$ structure, and use our obstruction theory to prove that a large class of spectra can be given an $A_\infty$ structure.

In section \ref{sec:cyclic} we define two notions of what it means for an $A_n$ structure on a spectrum $A$ to be ``commutative''. Let $K_n$ be the $n$'th associahedron, and let $W_n$ be the $n$'th cyclohedron. Then the $A_n$ structure is defined in terms of maps $(K_m)_+ \sma A^{(m)} \sar A$ for $m \leq n$. We say that the $A_n$ structure is \emph{cyclic} if the $n$ maps $(K_n)_+ \sma A^{(n)} \sar A$ obtained by first cyclically permuting the factors and then using the $A_n$ structure are homotopic, in the sense that a certain map $(\partial W_n)_+ \sma A^{(n)} \sar A$ which is given by the $A_n$ structure on each of the $n$ copies of $K_n$ on $\partial W_n$ can be extended to all of $W_n$. There is a dual notion, where instead of cyclically permuting the copies of $A$ we use the maps $(a_1,a_2,\ldots,a_n) \mapsto (a_2,\ldots,a_i,a_1,a_{i+1},\ldots,a_n)$. The notion of a cyclic $A_n$ structure plays a natural role in the study of topological Hochschild homology while the notion of a cocyclic $A_n$ structure plays a corresponding role for topological Hochschild cohomology.

In section $5$ we connect the ``commutativity'' of the multiplication on $A$ with $THH(A)$, and do some sample calculations. Baker and Lazarev (\cite[Theorem 3.1]{BaLa}) proved that
\begin{equation}
THH_{KU}(KU/2) \simeq KU^\wedge_2,
\end{equation}
and we are able to vastly generalize their result. For example, for an odd prime $p$, $THH_{KU}(KU/p)$ is not constant over the moduli space of $A_\infty$ structures on $KU/p$. For $p-1$ of the $p$ possible homotopy classes of multiplications ($A_2$ structures) on $KU/p$ we find that $THH_{KU}(KU/p) \simeq KU^\wedge_p$, while $THH_{KU}(KU/p)$ is a finite extension of $KU^\wedge_p$, ramified at $p$, on the rest of the moduli space. But the extension has degree at most $p-1$, so while $KU^\wedge_p \sar THH_{KU}(KU/p)$ might be a ramified extension, it is always tamely ramified.

In the last section, which is somewhat different from the rest of the paper, we compare $THH$ of Morava $K$-theory over different ground rings. While the calculations before used the corresponding Morava $E$-theory or Johnson-Wilson spectrum as the ground ring, we are really interested in using the sphere spectrum as the ground ring. We prove that in this particular case the choice of ground ring does not matter, that in fact the canonical maps
\begin{equation}
THH_R(A) \lar THH_S(A)
\end{equation}
and
\begin{equation}
THH^S(A) \lar THH^R(A)
\end{equation}
are weak equivalences.

\subsection*{Notation}
Throughout the paper $R$ will denote a commutative $S$-algebra as in \cite{EKMM}, and all smash products and function spectra will be over $R$ unless indicated otherwise. We will often assume that $R$ is even, meaning that $R_*$ is concentrated in even degrees. We let $I$ be a regular ideal in $R_*$, generated by the regular sequence $(x_1,\ldots,x_m)$, with $|x_i|=d_i$. For some applications we allow $I$ to be infinitely generated. We let $A$ be an $R$-module, often an $A_n$ algebra for some $2 \leq n \leq \infty$, and we let $M$ be another $R$-module, often an $(A_n,A)$-bimodule for some $2 \leq n \leq \infty$. We will often take $A=R/I$, an $R$-module with $A_*=R_*/I$.

We will be especially interested in writing the Morava $K$-theory spectra as quotients $R/I$ for suitable $R$ and $I$. To this end, let $\widehat{E(n)}$ denote the $K(n)$-localization of the Johnson-Wilson spectrum $E(n)$. It has homotopy groups
\begin{equation}
\widehat{E(n)}_* \cong \Z_{(p)}[v_1,\ldots,v_{n-1},v_n,v_n^{-1}]_I^\wedge,
\end{equation}
the $I$-completion of $E(n)_*$. Here $I=(p,v_1,\ldots,v_{n-1})$, and $K(n)=\widehat{E(n)}/I$ has homotopy groups
\begin{equation}
K(n)_* \cong \F_p[v_n,v_n^{-1}],
\end{equation}
which is a graded field with $|v_n|=2(p^n-1)$. It is known \cite{RoWh02} that $\widehat{E(n)}$ is an $E_\infty$ ring spectrum, or equivalently (\cite[Corollary II.3.6]{EKMM}) a commutative $S$-algebra.

We will also consider a $2$-periodic versions of Morava $K$-theory. Let $E_n=E_{(k,\Gamma)}$ be the Morava $E$-theory spectrum associated to a formal group $\Gamma$ of height $n$ over a perfect field $k$ of characteristic $p$. Then
\begin{equation}
(E_n)_* \cong \mathbb{W}k[[u_1,\ldots,u_{n-1}]][u,u^{-1}],
\end{equation}
where $\mathbb{W}k$ denotes the Witt vectors on $k$, to be thought of as $k[[u_0]]$. Here $|u_i|=0$, and we take $|u|=2$ rather than $-2$ as some authors do. Thus $u$ corresponds to the Bott element in the complex $K$-theory spectrum, rather than its inverse. It is known \cite{GoHo} that $E_n$ is an $E_\infty$ ring spectrum, or equivalently a commutative $S$-algebra. We let $I=(p,u_1,\ldots,u_{n-1})$ and $K_n=E_n/I$. Thus
\begin{equation}
(K_n)_* \cong k[u,u^{-1}].
\end{equation}

\subsection*{Acknowledgements}
An earlier version of this paper formed parts of the author's PhD thesis at the Massachusetts Institute of Technology under the supervision of Haynes Miller.

\section{$THH$ of $A_\infty$ ring spectra} \label{sec:THHdef}
We first recall some things from \cite{An_cyclic}. Let $\bD C$ be Connes' category of cyclic sets, and recall the construction in \cite[Definition 3.1]{An_cyclic} of the category $\bD C_\K$, which is a version of the category of cyclic sets which is enriched over $\Top$, the category of (compactly generated, weak Hausdorff) topological spaces. Here $\K$ denotes the associahedra operad (by operad we really mean non-$\Sigma$ operad), with $\K(\{1,2,\ldots,n\})=K_n \cong D^{n-2}$, and each $Hom$ space in $\bD C_\K$ is a disjoint union of products of associahedra. Let $^0 \! \bD C$ be the category of based cyclic sets and basepoint-preserving maps, which is equivalent to the simplicial indexing category $\bD^{op}$ \cite[Lemma 3.3]{An_cyclic}, and let $^0 \! \bD C_\K$ be the corresponding enriched category.

Let $\mathcal{W}$ be the collection of cyclohedra, with $\mathcal{W}(\{0,1,\ldots,n-1\})=W_n \cong D^{n-1}$. Here we think of $\{0,1,\ldots,n-1\}$ as a cyclically ordered set, i.e., an object in $\bD C$. The cyclohedra $\mathcal{W}$ assemble to a functor $\bD C_\K^{op} \sar \Top$, though we will only need that $\mathcal{W}$ is a functor $^0 \! \bD C_\K^{op} \sar \Top$ for defining $THH$.

Let $\M_S$ be the category of $S$-modules as in \cite{EKMM}. By a spectrum we mean an $S$-module, and we write $\sma$ for $\sma_S$. Suppose that $A \in \M_S$ is an $A_\infty$ ring spectrum, by which we mean an algebra over the operad $\K$. Thus $A$ comes with maps
\begin{equation}
\phi_n : (K_n)_+ \sma A^{(n)} \sar A
\end{equation}
for $n \geq 0$ making certain diagrams commute. Also suppose that $M$ is an $A$-bimodule, a.k.a.~an $(A_\infty,A)$-bimodule. By that we mean that there are maps \begin{equation}
\xi_{n,i} : (K_n)_+ \sma A^{(i-1)} \sma M \sma A^{(n-i)} \sar M
\end{equation}
for $n \geq 1$ and $1 \leq i \leq n$ making similar diagrams commute. Here, and throughout the paper, $A^{(j)}$ denotes the $j$-fold smash product of $A$ with itself.

In this situation we can define a functor $B^{cy}(A;M) : {}^0 \! \bD C_\K \sar \M_S$ by $B^{cy}(A;M)(\{0,1,\ldots,n\})=M \sma A^{(n)}$, and a functor $C^{cy}(A;M) : {}^0 \! \bD C_\K^{op} \sar \M_S$ by $C^{cy}(A;M)(\{0,1,\ldots,n\})=F_S(A^{(n)},M)$.

\begin{definition}
Let $A$ be an $A_\infty$ ring spectrum and let $M$ be an $A$-bimodule, with $A$ and $M$ cofibrant (q-cofibrant in the terminology of \cite{EKMM}) in $\M_S$. Topological Hochschild homology of $A$ with coefficients in $M$ is the spectrum
\begin{equation}
THH^S(A)=\mathcal{W} \otimes_{{}^0 \! \bD C_\K} B^{cy}(A;M)
\end{equation}
while topological Hochschild cohomology of $A$ with coefficients in $M$ is the spectrum
\begin{equation}
THH_S(A)=Hom_{{}^0 \! \bD C_\K}(\mathcal{W},C^{cy}(A;M)).
\end{equation}
\end{definition}

Here $-\otimes_{{}^0 \! \bD C_\K}\!\!-$ and $Hom_{{}^0 \! \bD C_\K}(-,-)$ are defined as a suitable coequalizer and equalizer, see \cite[Definition 3.7]{An_cyclic}.

\begin{remark}
If $A=M$ then $B^{cy}(A;A)$ can be extended to a functor from $\bD C_\K$ to $\M_S$. Because $\mathcal{W}$ is a functor $\bD C_\K^{op} \sar \Top$, one can show that $THH^S(A)$ has an action of $S^1$, in much the same way as in the classical situation. We omit the details, since we will not need the $S^1$-action.
\end{remark}

If $R$ is a commutative $S$-algebra, $A$ is an $A_\infty$ $R$-algebra and $M$ is an $A$-bimodule, we can define $THH^R(A;M)$ and $THH_R(A;M)$ by taking all smash products and function spectra in the category $\M_R$ of $R$-modules instead of $\M_S$.

To avoid over-using the word cofibrant, we will assume that all our spectra are cofibrant in $\M_R$.

If we have a map $M \sar M'$ of $A$-bimodules, we get maps $THH^R(A;M) \sar THH^R(A;M')$ and $THH_R(A;M) \sar THH_R(A;M')$. If we have a map $A \sar A'$ of $A_\infty$ ring spectra and $M'$ is an $A'$-bimodule, we get maps $THH^R(A;M') \sar THH^R(A';M')$ and $THH_R(A';M') \sar THH_R(A;M')$.

\begin{proposition}
If $A$ is strictly associative, our definition of $THH^R(A;M)$ agrees (in $\M_R$) with the definition of $thh^R(A;M)$ given in \cite[IX.2]{EKMM}. (They only define $THH_R(A;M)$ in the derived category.) Moreover, our definition is homotopy invariant in the following sense: If $(A,M) \sar (A',M')$ is a map of $A_\infty$ ring spectra and bimodules such that $A \sar A'$ and $M \sar M'$ are weak equivalences, then we have weak equivalences
\begin{equation}
THH^R(A;M) \arrow{\simeq} THH^R(A;M') \arrow{\simeq} THH^R(A',M')
\end{equation}
and
\begin{equation}
THH_R(A;M) \arrow{\simeq} THH_R(A;M') \overset{\simeq}{\longleftarrow} THH_R(A';M').
\end{equation}
\end{proposition}

\begin{proof}
The first claim follows in the same way as \cite[Proposition 4.13]{An_cyclic}, and the homotopy invariance follows from the theory of enriched Reedy categories developed in \cite{An_R}.
\end{proof}

\begin{remark}
If $A$ and $M$ are not cofibrant, the correct way to define $THH^R(A;M)$ is as the ``geometric realization'' of a cofibrant replacement of $B^{cy}(A;M)$ in the Reedy model category on functors $^0 \! \bD C_\K \sar \M_R$, and the correct way to define $THH_R(A;M)$ is as ``Tot'' of a fibrant replacement of $C^{cy}(A;M)$.
\end{remark}

For ease of reference, we recall the standard spectral sequences used to calculate the homotopy or homology groups of $THH$.

\begin{proposition} \label{canonicalSS} (\cite[chapter IX]{EKMM})
There are spectral sequences
\begin{eqnarray}
E^2_{s,t}=Tor^{\pi_*(A \sma_R A^{op})}_{s,t}(A_*,M_*) \Longrightarrow \pi_{s+t} THH^R(A;M), \label{homologySS} \\
E_2^{s,t}=Ext_{\pi_*(A \sma_R A^{op})}^{s,t}(A_*,M_*) \Longrightarrow \pi_{t-s} THH_R(A;M). \label{cohomologySS}
\end{eqnarray}
If $E$ is a commutative $R$-algebra, or if $E_*(A \sma_R A^{op})$ is flat over $\pi_*(A \sma_R A^{op})$, there are spectral sequences
\begin{eqnarray}
E^2_{s,t}=Tor^{E_*(A \sma_R A^{op})}_{s,t}(E_*^RA,E_*^RM) \Longrightarrow E_{s+t} THH^R(A;M), \\
E_2^{s,t}=Ext_{E_*(A \sma_R A^{op})}^{s,t}(E_*^RA,E_*^RM) \Longrightarrow E_{t-s} THH_R(A;M).
\end{eqnarray}
Here $E_*^RX$ means $\pi_*(E \sma_R X)$.
\end{proposition}

Under reasonable finiteness conditions on each group these spectral sequences converge strongly (\cite[Theorem 6.1 and 7.1]{Bo99}).

The spectral sequence 
\begin{equation}
E^2_{*,*}=Tor^{H_*(A \sma_S A^{op};\F_p)}_{*,*}(H_*(A;\F_p);H_*(M;\F_p)) \Longrightarrow H_*(THH^S(A;M);\F_p)
\end{equation}
is called the B\"okstedt spectral sequence, after Marcel B\"okstedt who first defined topological Hochschild homology \cite{Bo1}, \cite{Bo2}.

Topological Hochschild cohomology acts on topological Hochschild homology via maps
\begin{equation}
 THH_R(A) \sma_R THH^R(A;M) \lar THH^R(A;M).
\end{equation}
This is clear from the definition of $THH_R(A)$ and $THH^R(A;M)$ as the derived spectra $F_{A \sma_R A^{op}}(A,A)$ and $M \sma_{A \sma_R A^{op}} A$ from \cite[IX.1]{EKMM}. It is also clear that the natural action
\begin{equation}
 E_2^{s,t}(A) \otimes_{R_*} E^2_{p,q}(A;M) \lar E^2_{p-s,q+t}(A;M)
\end{equation}
given by the action of $Ext$ on $Tor$ converges to the action of $\pi_*THH_R(A)$ on $\pi_*THH^R(A;M)$.

If $A$ is only an $A_n$ ring spectrum and $M$ is an $(A_n,A)$-bimodule, meaning that we have coherent maps $(K_m)_+ \sma A^{(i-1)} \sma M \sma A^{(m-i)} \sar M$ for $m \leq n$ and $1 \leq i \leq m$, we get a functor, which we will still denote by $B^{cy}(A;M)$, from ${}^0 \! \bD C_{\K_n}$ to $\M_R$. Here $\K_n$ is the operad generated by $K_m$ for $m \leq n$. We can then consider the functor $\mathcal{W}_n : \bD C_{\K_n} \sar \Top$ generated by $W_m$ for $m \leq n$ and define spectra
\begin{equation}
sk_{n-1} THH^R(A;M)= \mathcal{W}_n \otimes_{{}^0 \! \bD C_{\K_n}} B^{cy}(A;M)
\end{equation}
and
\begin{equation}
Tot^{n-1} THH_R(A)=Hom_{{}^0 \! \bD C_{\K_n}}(\mathcal{W}_n,C^{cy}(A;M)).
\end{equation}

If the $A_n$ structure on $A$ can be extended to an $A_\infty$ structure, and the $(A_n,A)$-bimodule structure on $M$ can be extended to an $(A_\infty, A)$-bimodule structure, then these constructions coincide with the constructions obtained from the skeletal and total spectrum filtrations.

\section{$A_\infty$ obstruction theory} \label{sec:Ainftyobtheory}
In this section we set up an obstruction theory for endowing a spectrum $A$ with an $A_\infty$ structure. There are basically two ways to do this, both reducing to topological Hochschild cohomology calculations. One way, which only works for connective spectra, is to build an $A_\infty$ structure by induction on the Postnikov sections $P_m A$ \cite{DuSh}. The other, which we will use here, is to proceed by induction on the $A_n$ structure.

The original reference for this is \cite{Ro88}, but Robinson implicitly assumes that the multiplication on $A$ is homotopy commutative, an assumption we would very much like to get rid of. Other works on the subject, such as \cite{Go}, also assumes that the multiplication is homotopy commutative.

The space $Aut(A)$ acts on the set (or space) of $A_\infty$ structures by conjugation, and it is natural to consider two $A_\infty$ structures on $A$ equivalent if they differ by conjugation by an automorphism of $A$. We will not attempt to mod out by the action of $Aut(A)$ here (except in Remark \ref{rem:autos}), though we hope to come back to this elsewhere. Corollary \ref{cor:samespaceofAinfty} is certainly not true after modding out by $Aut(A)$, as $Aut_R(A)$ is very different from $Aut_S(A)$ in this case.

The results in this section strengthen several results already in the literature. For example, Corollary \ref{cor:R/xisainfty} strengthens \cite[Proposition 3.1(1)]{Str}, which says that for $R$ even commutative and $x$ a nonzero divisor, any multiplication on $R/x$ is homotopy associative, to saying that any multiplication on $R/x$ can be extended to an $A_\infty$ structure. Corollary \ref{cor:evenquotientainfty} strengthens various results in \cite{La03} and \cite{BaJe} about the associativity of $MU/I$ for certain regular ideals to all regular ideals.

The most important change from \cite{Ro88} is that instead of $\pi_* A \sma A$ we consider $\pi_* A \sma A^{op}$. The result is that the obstructions lie in the $(-3)$-stem of the potential spectral sequence converging to $\pi_* THH_S(A)$, in a sense we will make more precise below.

The obstruction theory works just as well in the category $\M_R$ for a commutative $S$-algebra $R$, in which case we are looking for $A_\infty$ $R$-algebra structures on $A$. As an example of the power of this obstruction theory, we prove that if $R$ is even and $I$ is a regular ideal in $R_*$, $R/I$ can always be given an $A_\infty$ $R$-algebra structure.

Suppose that we have an $A_{n-1}$ structure on a spectrum $A$ and we want to extend it to an $A_n$ structure. Then we need a map
\begin{equation}
(K_n)_+ \sma A^{(n)} \sar A
\end{equation}
which is compatible with the $A_{n-1}$ structure. Because all the faces of $K_n$ are products of associahedra of lower dimension, the map $(K_n)_+ \sma A^{(n)} \sar A$ is determined on $\partial K_n \sma A^{(n)} \simeq \Sigma^{n-3} A^{(n)}$. Thus the obstruction to extending the given $A_{n-1}$ structure to an $A_n$ structure lies in
\begin{equation}
[\Sigma^{n-3} A^{(n)},A] = A^{3-n}(A^{(n)}).
\end{equation}

The unitality condition on the $A_n$ structure also fixes the map on $(K_n)_+ \sma s_j A^{(n-1)}$ for $0 \leq j \leq n-1$, where $s_j : A^{(n-1)} \sar A^{(n)}$ is given by the unit $R \sar A$ on the appropriate factor. (This does not quite make sense; what we mean is that the appropriate diagram is required to commute.) If we define $\bar{A}$ as the cofiber of the unit map $R \sar A$, we can then say that the obstruction lies in
\begin{equation}
[\Sigma^{n-3} \bar{A}^{(n)},A]=A^{3-n}(\bar{A}^{(n)}).
\end{equation}

We also note that if the set of homotopy classes of $A_n$ structures on $A$ with a fixed $A_{n-1}$ structure is nonempty, it is isomorphic to $A^{2-n}(\bar{A}^{(n)})$ as a set. This set has no group structure, but we can say that it is an $A^{2-n}(\bar{A}^{(n)})$-torsor.

We define a bigraded group $E_1^{*,*}$ by
\begin{equation}
E_1^{s,t}=A^{-t}(A^{(s)}).
\end{equation}
To take the unitality condition into account we also define
\begin{equation}
\bar{E}_1^{s,t}=A^{-t}(\bar{A}^{(s)}).
\end{equation}
Thus the obstruction to extending a given $A_{n-1}$ structure to an $A_n$ structure lies in $\bar{E}_1^{n,n-3}$.

Note that we do not need the existence of a homotopy associative multiplication on $A$, which is the data Robinson starts with in \cite{Ro88}, to define $\bar{E}_1^{*,*}$. We get the following:

\begin{theorem} \label{AinftyobE_1}
Given an $A_{n-1}$ structure on $A$, $n \geq 2$, the obstruction to the existence of an $A_n$ structure extending the given $A_{n-1}$ structure lies in $\bar{E}_1^{n,n-3}$.
\end{theorem}

\begin{proof}
This is clear from the above discussion when $n \geq 3$, we give a separate argument for $n=2$.

There is a cofiber sequence $\Sigma^{-1} \bar{A} \sma \bar{A} \sar A \vee A \sar A \sma A$. Consider the fold map $A \vee A \sar A$. It can be extended to a map $A \sma A \sar A$ if and only if the composite $\Sigma^{-1} \bar{A} \sma \bar{A} \sar A \vee A \sar A$ is null, so the obstruction lies in $[\Sigma^{-1} \bar{A} \sma \bar{A},A]=\bar{E}_1^{2,-1}$. If the fold map can be extended to $A \sma A$, it will automatically be unital, so it is an $A_2$ structure.
\end{proof}

Let $R$ be an even commutative $S$-algebra. Because $R$ is not a cell $R$-module, we use the sphere $R$-modules $S^n_R$ to make cell $R$-modules, as in \cite[III.2]{EKMM}. Given $x \in \pi_d R$, we define $R/x$ as the cofiber
\begin{equation}
S^d_R \overset{x}{\sar} S^0_R \sar R/x.
\end{equation}

\begin{corollary} \label{cor:R/xisainfty}
Let $R$ be an even commutative $S$-algebra, and let $x \in \pi_d R$ be a nonzero divisor. Then any $A_{n-1}$ structure on $A=R/x$ can be extended to an $A_n$ structure, for any $n \geq 2$. In particular, the set of $A_\infty$ structures on $A$ is nonempty.
\end{corollary}

\begin{proof}
In this case $\overline{R/x} \cong \Sigma^{d+1} R$. Thus 
\begin{equation}
\bar{E}_1^{s,t}=R^{-t}((\Sigma^{d+1} R)^{(s)}) \cong R^{-t}(\Sigma^{s(d+1)} R) \cong \pi_{s(d+1)+t} R,
\end{equation}
and the obstruction lies in $\bar{E}_1^{n,n-3}=\pi_{n(d+1)+n-3} R$, which is zero because $R$ (and $d$) is even.
\end{proof}

In particular, this settles \cite[Conjecture 2.16]{BaLa}, where Baker and Lazarev conjecture that any homotopy associative multiplication on $R/x$ can be extended to an $A_\infty$ multiplication.

\begin{example}
As an example of how $R/x$ fails to be $A_\infty$ when $R$ is not even, let us consider the case $R=S$ and $x=p$, so $S/p=M_p$ is the mod $p$ Moore spectrum. In this case the obstruction to an $A_n$ structure lies in $\pi_{2n-3} M_p$, which is zero for $n<p$, but $\pi_{2p-3} M_p \cong \Z/p$ is generated by $S^{2p-3} \overset{\alpha_1}{\sar} S^0 \sar M_p$.

The obstruction is in fact nonzero. One way to show this is to consider the map $M_p \sar H\Z/p$. If $M_p$ is $A_p$, then this is a map of $A_p$ ring spectra, and the induced map $H_*(M_p,\Z/p) \sar H_*(H \Z/p;\Z/p)$ commutes with $p$-fold Massey products. But there is a $p$-fold Massey product $\langle \bar{\tau}_i,\ldots,\bar{\tau}_i \rangle=-\bar{\xi}_{i+1}$ in $H_*(H\Z/p;\Z/p)=A_*$ defined with no indeterminacy,\footnote{We have not found this statement in the literature, but the proof is easy. By Kochman \cite[Corollary 20]{Ko72}, the $p$-fold Massey product on a class $x$ in dimension $2n-1$ is given by $-\beta Q^n(x)$, and by Steinberger's calculations \cite[Theorem III.2.3]{BMMS} this gives the result.} and in particular the image of the generator $a \in H_1(M_p;\Z/p)$ supports a nonzero $p$-fold Massey product while $a$ clearly does not.
\end{example}

Now suppose that $A$ comes with an $A_3$ structure, i.e., a unital map $\phi_2 : A \sma A \sar A$ and a homotopy $\phi_3$ from $\phi_2(\phi_2 \sma 1)$ to $\phi_2(1 \sma \phi_2)$. Then there are $s+2$ maps $A^t(A^{(s)}) \sar A^t(A^{(s+1)})$, which we denote by $d^i$ for $0 \leq i \leq s+1$. Here $d^0$ sends $f : A^{(s)} \sar A$ to $A^{(s+1)} \overset{1 \sma f}{\sar} A^{(2)} \overset{\phi_2}{\sar} A$, $d^i$ sends $f$ to $A^{(s+1)} \arrow{1^{i-1} \sma \phi_2 \sma 1^{s-i}} A^{(s)} \overset{f}{\sar} A$ for $1 \leq i \leq s$ and $d^{s+1}$ sends $f$ to $A^{(s+1)} \overset{f \sma 1}{\sar} A^{(2)} \overset{\phi_2}{\sar} A$.

Adding the obvious codegeneracy maps, this structure makes $E_1^{*,*}$ into a graded cosimplicial group. Note that we could not do this with only an $A_2$ structure, because homotopy associativity is needed to make sure the cosimplicial identities hold.

With this construction, $\bar{E}_1^{*,*}$ is the associated normalized cochain complex, with differential $d=\sum (-1)^i d^i$. We let $E_2^{*,*}$ be the homology of $(\bar{E}_1^{*,*},d)$. The existence of an $A_3$ structure is also needed to make sure $d^2=0$.

\begin{lemma}
Suppose we have an $A_{n-1}$ structure on $A$, $n \geq 4$. Let $c_n$ be the obstruction to the existence of an $A_n$ structure extending the given $A_{n-1}$ structure. Then $d(c_n)=0$.
\end{lemma}

\begin{proof}
The obstruction $c_n$ is a map $\partial K_n \sma A^{(n)} \sar A$, so we think of $c_n$ as the boundary of $K_n$. In this way we can think of $d^i(c_n)$ as the boundary of one of the copies of $K_n$ on $\partial K_{n+1}$, which is a codimension $2$ subcomplex of $K_{n+1}$, and we can consider $d(c_n)$ as a formal sum of codimension $2$ subcomplexes of $K_{n+1}$. The faces that lie in the intersection of two copies of $K_n$ sum to zero, while the rest are null because we can fill the copies of $K_i \times K_{n-i+2}$ on $\partial K_{n+1}$ for $3 \leq i \leq n-1$.
\end{proof}

If we change the $A_{n-1}$ structure by $f$, the obstruction to an $A_n$ structure changes by $df$. Thus we get the following:

\begin{theorem} \label{AinftyobE_2} (Compare \cite[Theorem 1.11]{Ro88})
Suppose we have an $A_{n-1}$ structure on $A$, $n \geq 4$. The obstruction to the existence of an $A_n$ structure on $A$, while allowing the $A_{n-1}$ structure to vary but fixing the $A_{n-2}$ structure lies in $E_2^{n,n-3}$.
\end{theorem}

Under very reasonable conditions on $A$, for example if $\pi_* A \sma_R A$ is projective over $A_*$, we can identify $E_2^{*,*}$ with $Ext_{\pi_*A \sma_R A^{op}}^{*,*}(A_*,A_*)$. Now we are in a position to prove that any homotopy associative multiplication on $A=R/I$, where $R$ is even and $I$ is a regular ideal, can be extended to an $A_\infty$ structure. Here $R/I$ is defined as follows. Let $I=(x_1,x_2,\ldots)$, where $(x_1,x_2,\ldots)$ is a regular sequence. Then $R/I$ is the (possibly infinite) smash product $R/x_1 \sma_R R/x_2 \sma_R \ldots$. It is clear (\cite[Corollary V.2.10]{EKMM}) that $R/I$ does not depend on the choice of regular sequence generating $I$. First we need to know the structure of $\pi_* A \sma_R A^{op}$. Let $d_i$ be the degree of $x_i$.

\begin{proposition} \label{prop:AAopring}
Given any homotopy associative multiplication on $A=R/I$ with $R$ even and $I=(x_1,x_2,\ldots)$ a regular ideal, $\pi_* A \sma_R A^{op}$ is given by
\begin{equation}
\pi_*A \sma_R A^{op}=\Lambda_{A_*}(\alpha_1,\alpha_2,\ldots)
\end{equation}
as a ring. Here $|\alpha_i|=d_i+1$.
\end{proposition}

\begin{proof}
This is well known (\cite{BaJe}, \cite{La03}). The proofs in \cite{BaJe} and \cite{La03} both use that $\pi_*F_R(A,A)$ is a (completed) exterior algebra together with a Kronecker pairing. Here we present a different proof:

There is a multiplicative K\"unneth spectral sequence (see \cite{BaLa01})
\begin{equation}
E^2_{*,*}=Tor^{R_*}_{*,*}(A_*, A^{op}_*) \Longrightarrow \pi_* A \sma_R A^{op}.
\end{equation}
By using a Koszul resolution of $A_*=R_*/I$ it is easy to see that $E^2_{*,*}=\Lambda_{A_*}(\alpha_1,\alpha_2,\ldots)$ with $\alpha_i$ in bidegree $(1,d_i)$. The spectral sequence collapses, so all we have to do is to show that there are no multiplicative extensions. Because $\alpha_i^2$ is well defined up to \emph{lower} filtration and $E^2_{1,*}$ is concentrated in odd total degree, it follows that $\alpha_i^2 \in A_* \otimes_{R_*} A_*^{op} \cong A_*$ in $\pi_*A \sma_R A^{op}=A^R_*A^{op}$. Now there are several ways to show that $\alpha_i^2=0$. If we denote the map $A^R_*A^{op} \lar A_*$ by $\epsilon$, it is enough to show that $\epsilon(\alpha^2)=0$ since $\epsilon$ gives an isomorphism from filtration $0$ in the spectral sequence to $A_*$. For example, we can use that $A$ is an $A \sma_R A^{op}$-module and study the two maps $A^R_*A^{op} \otimes A^R_*A^{op} \otimes A_* \lar A_*$. One sends $\alpha_i \otimes \alpha_i \otimes 1$ to $\epsilon(\alpha_i^2)$, the other one sends it to $0$.
\end{proof}

An extension of the argument in the proof shows that there cannot even be any Massey products in $\pi_*A \sma_R A^{op}$, by comparing brackets formed in $(A \sma_R A^{op})^{(n)}$ and $(A \sma_R A^{op})^{(n-1)} \sma A$.

The above result is \emph{not} true for $A \sma_R A$, in which case $\alpha_i$ might very well square to something non-zero.

\begin{corollary} \label{cor:evenquotientainfty}
Suppose $A=R/I$ with $R$ even and $I$ regular has an $A_{n-1}$ structure, $n \geq 4$. Then $A$ has an $A_n$ structure with the same underlying $A_{n-2}$ structure. In particular, any homotopy associative multiplication on $A=/I$ can be extended to an $A_\infty$ structure.
\end{corollary}

\begin{proof}
Using Theorem \ref{AinftyobE_2}, the relevant obstructions lie in $Ext_{\pi_* A \sma_R A^{op}}^{n,n-3}(A_*,A_*)$. In particular the obstructions are in odd total degree. But $\pi_* A \sma_R A^{op} \cong \Lambda_{A_*}(\alpha_1,\alpha_2,\ldots)$ with $|\alpha_i|=|d_i|+1$, so $Ext$ over it is a polynomial algebra 
\begin{equation}
Ext_{\pi_* A \sma_R A^{op}}^{*,*}(A_*,A_*) \cong A_*[q_1,q_2,\ldots] \label{eq:THHR/ISS}
\end{equation}
with $|q_i|=(1,-d_i-1)$. Thus $Ext_{\pi_* A \sma_R A^{op}}^{*,*}(A_*,A_*)$ is concentrated in even total degree, and there can be no obstructions.
\end{proof}

Equation \ref{eq:THHR/ISS} in the above proof also gives the $E_2$-term of the canonical spectral sequence calculating $\pi_* THH_R(R/I)$. A similar calculation gives the $E^2$-term of the spectral sequence calculating $\pi_* THH^R(R/I)$ as a divided power algebra, as in equation \ref{eq:THCSS} and \ref{eq:THHSS}.
\begin{remark}

It might seem like Corollary \ref{cor:evenquotientainfty} follows from Corollary \ref{cor:R/xisainfty}, because an $A_\infty$ structure on each $R/x_i$ gives an $A_\infty$ structure on $R/I$, but there are multiplications on $R/I$ which do not come from smashing together multiplications on each $R/x_i$, so this is a stronger result.

There might also be $A_2$ structures on $R/I$ which do not extend to $A_3$, although any $A_2$ structure obtained by smashing together $A_2$ structures on each $R/x_i$ will certainly be homotopy associative.
\end{remark}

We will need a more precise classification of the $A_2$ structures on $R/I$ which can be extended to $A_3$, and hence to $A_\infty$. Recall (\cite[Proposition 4.15]{Str}) that $A_R^*A$ is a (completed) exterior algebra
\begin{equation}
A_R^*A \cong \hat{\Lambda}_{A_*}(Q_1,Q_2,\ldots),
\end{equation}
where $Q_i$ is obtained from the composite $R/x_i \overset{\beta_i}{\sar} \Sigma^{d_i+1} R \sar \Sigma^{d_i+1} R/x_i$ and has degree $-d_i-1$.

\begin{theorem}
Fix a homotopy associative multiplication $\phi^0$ on $A=R/I$. Given any other homotopy associative multiplication $\phi$ on $A$, it can be written uniquely as
\begin{equation} \label{allass}
\phi=\phi^0 \prod_{i,j} \big(1 \sma 1+v_{ij} Q_i \sma Q_j \big)
\end{equation}
for some $v_{ij} \in \pi_{d_i+d_j+2} A$, where the product denotes composition (which can be taken in any order, because all the factors are even). Conversely, any $\phi$ that can be written in this form is homotopy associative.
\end{theorem}

\begin{proof}
Associativity is some kind of cocycle condition, and one could imagine a simple proof based on this. However, the relevant maps $A^0(A \sma A) \lar A^0(A \sma A \sma A)$ are not linear, and this complicates things.

We use the K\"unneth isomorphism
\begin{equation}
A^*A \cong Hom_{A_*}(A_*A,A_*)
\end{equation}
and similar formulas for $A^*(A^{(2)})$ and $A^*(A^{(3)})$. These isomorphisms depend on a choice of multiplication, and we will use $\phi^0$ for each of them. For example, the map $A^*A \lar Hom_{A_*}(A_*A,A_*)$ is given by sending $A \arrow{f} A$ to $A_*A \arrow{A_*f} A_*A \arrow{\phi^0} A_*$.

Let $\epsilon : A_* A \lar A_*$ be the map induced by $\phi^0$. To check if $\phi$ is associative, it is enough to check whether or not the diagram
\begin{equation} \label{checkass}
\xymatrix{
A_*A \otimes_{A_*} A_*A \otimes_{A_*} A_*A \ar[rr]^-{\phi \sma 1} \ar[dd]_{1 \sma \phi} & & A_*A \otimes_{A_*} A_*A \ar[d]^\phi \\
 & & A_*A \ar[d]^\epsilon \\
A_*A \otimes_{A_*} A_*A \ar[r]^-\phi & A_*A \ar[r]^-\epsilon & A_*}
\end{equation}
commutes.

Recall that $A^*A \cong \hat{\Lambda}_{A_*}(Q_1,\ldots,Q_n)$ and that $A_*A \cong \Lambda_{A_*}(\alpha_1,\ldots,\alpha_n)$, at least additively. Under the K\"unneth isomorphism $Q_i$ corresponds to the map sending $\alpha_i$ to $1$.

Now suppose that $\phi$ is some unital product on $A$. We can write
\begin{eqnarray}
\phi=\phi^0 \prod_{I,J} \big( 1 \sma 1+v_{IJ} Q_I \sma Q_J \big),
\end{eqnarray}
where $I$ and $J$ run over indexes $I=(i_1,\ldots,i_r)$ and $J=(j_1,\ldots,j_s)$, where $Q_I=Q_{i_1} \cdots Q_{i_r}$ and $Q_J=Q_{j_1} \cdots Q_{j_s}$. Let $|I|$ denote the number of indices in $I$. By unitality we have $|I| > 0$ and $|J| > 0$, and because $A_*$ is even $|I|+|J|$ has to be even.

If $\phi=\phi^0(1 \sma 1+v_{ij} Q_i \sma Q_j)$, then we can calculate $\phi(\phi \sma 1)$ and $\phi(1 \sma \phi)$ using diagram \ref{checkass}. For example, $\phi(\phi \sma 1)$ and $\phi(1 \sma \phi)$ both send $\alpha_i \otimes \alpha_j \otimes 1$ to $v_{ij}$, as we see by following diagram \ref{checkass} around both ways. Similarly, they send $\alpha_i \otimes 1 \otimes \alpha_j$ and $1 \otimes \alpha_i \otimes \alpha_j$ to $v_{ij}$, and they send $\alpha_i \otimes \alpha_i \alpha_j \otimes \alpha_j$ to $-v_{ij}^2$. Those are all the relevant terms, and shows that
\begin{multline}
\phi(1 \sma \phi)=\phi(\phi \sma 1)=\phi^0(\phi^0 \sma 1) \circ \\
\Big( v_{ij}(Q_i \sma Q_j \sma 1+Q_i \sma 1 \sma Q_j+1 \sma Q_i \sma Q_j)-v_{ij} Q_i \sma Q_{ij} \sma Q_j \Big).
\end{multline}
This shows that any $\phi$ as in the theorem is associative.

To show that none of the other products are associative, it is enough to show that
\begin{equation}
\phi=\phi^0(1 \sma 1+v_{IJ} Q_I \sma Q_J)
\end{equation}
is not associative for any $I$, $J$ with $|I|+|J|>2$. For example, if
\begin{equation}
\phi=\phi^0(1 \sma 1+ vQ_{ij} \sma Q_{kl})
\end{equation}
then $\phi(1 \sma \phi)$ sends $\alpha_i \alpha_j \otimes \alpha_k \otimes \alpha_l$ to $v$ but $\phi(\phi \sma 1)$ sends it to zero.
\end{proof}

\begin{remark} \label{unique_form_remark}
Alternatively, we can say that given a homotopy associative multiplication $\phi$ on $A$, it can be written as
\begin{equation}
\phi=\phi^0 \prod_{i \neq j} \big(1 \sma 1+v_{ij} Q_i \sma Q_j \big)
\end{equation}
for a unique $\phi^0$ which is obtained by smashing together multiplications on each $R/x_i$.
\end{remark}

By allowing the $A_{n-1}$ structure but fixing the $A_{n-2}$ structure, we have seen that the obstruction $c_n$ to an $A_n$ structure lies in $E_2^{n,n-3}$. In fact, we can do even better. By allowing the $A_{n-i}$ structure to vary for $1 \leq i \leq r-1$, the obstruction to the existence of an $A_n$ structure actually lies in $E_r^{n,n-3}$, provided that $n \geq 2r$.

The reason for this restriction on $n$ and $r$ is the following. The definition of $d_{r-1}(c_n)$ uses the $A_i$ structure for $i \leq r$, and if we change the $A_r$ structure, we change the definition of $d_{r-1} : E_{r-1}^{n,n-3} \sar E_{r-1}^{n+r-1,n+r-3}$. This nonlinear behavior prevents $d_{r-1}$ from squaring to zero, so to define $E_r^{*,*}$ we have to fix the $A_r$ structure.

\begin{theorem} \label{AinftyobE_r}
Suppose $n \geq 2r$. Then the obstruction to an $A_n$ structure on $A$, while allowing the $A_{n-r+1}$ through $A_{n-1}$ structure to vary but fixing the $A_{n-r}$ structure, lies in $E_r^{n,n-3}$. 
\end{theorem}

\begin{proof}
With $n \geq 2r$, this works just like Theorem \ref{AinftyobE_1} and \ref{AinftyobE_2}, except the definition of $d_i$ for $i \geq 2$ is slightly more complicated. Let $c_n$ be the obstruction to an $A_n$ structure. Then, if $f \in E_2^{n-2,n-4}$, $d_2(f)$ changes the obstruction by a sum of two terms. First, the obstruction problem for the $A_{n-1}$ structure changes, so the obstruction changes on the faces of $K_n$ of the form $K_{n-1}$. Second, the obstruction changes on the faces of $K_n$ of the form $K_3 \times K_{n-2}$ because the map depends directly on the $A_{n-2}$ structure on those faces. The general case is similar.

Proving that $d_i(c_n)=0$ for $i \geq 2$ is also a bit more complicated. We have already seen that $d_1(c_n)=0$. If we think of this as a map from a subcomplex of $K_{n+1}$, we can extend this to all of $\partial K_{n+1}$ without changing the map on faces of the form $K_i \times K_{n+2-i}$ for $3 \leq i \leq n-1$. This gives us a map from $\partial K_{n+1}$. Now $d_2(c_n)$ is given by a sum of two terms. First, it changes on the codimension $2$ faces of $K_{n+2}$ of the form $\partial K_{n+1}$ given by $d_1$ of the map from $\partial K_{n+1}$ we just found. Second, it changes on the faces of the form $K_3 \times \partial K_n$.

Now some parts cancel, and the rest are null because we can fill $K_i \times K_{n+3-i}$ for $4 \leq i \leq n-1$. The general case is similar.
\end{proof}

It is also possible to set up an obstruction theory for extending a map $f: A \sar B$ between $A_\infty$ ring spectra to an $A_\infty$ map. We give a brief outline of one way to do this. We need to specify what we mean by an $A_\infty$ map. Requiring the diagrams
\[ \xymatrix{
(K_n)_+ \sma A^{(n)} \ar[r] \ar[d] & (K_n)_+ \sma B^{(n)} \ar[d] \\
A \ar[r] & B} \]
to commute on the nose is too restrictive.

Instead we should require, for example, that there is a homotopy between the two ways of going from $A \sma A$ to $B$, and then higher homotopies between all the ways of going from $A^{(n)}$ to $B$. This is encoded in what is sometimes called a colored operad, or a many-sorted operad, or a multicategory \cite{Leinster}. The idea goes back to Boardman and Vogt \cite{BoVo73}, who considered colored PROPs. In this case we want a multicategory $\K_{0 \sar 1}$ with two objects satisfying the following conditions. First of all, each space where all the source objects and the target object agree gives an associahedron, i.e., $\K_{0 \sar 1}(\epsilon,\ldots,\epsilon;\epsilon)=K_n$ for $\epsilon=0,1$. Second of all, $\K_{0 \sar 1} (\epsilon_1,\ldots,\epsilon_n;0)$ is empty if some $\epsilon_i=1$. And finally, $\K_{0 \sar 1}(0,0;1)=I$ is an interval and each $\K_{0 \sar 1}(0,\ldots,0;1)$ is contractible.

An algebra over $\K_{0 \sar 1}$ is precisely a pair of $A_\infty$ algebras and an $A_\infty$ map between them. It is a little bit harder to describe the spaces in $\K_{0 \sar 1}$, but it is true that $\K_{0 \sar 1}(0,\ldots,0;1)$ (with $n$ inputs) is homeomorphic to $D^{n-1}$. Then the usual obstruction theory argument shows that the obstruction to extending an $A_{n-1}$ map from $A$ to $B$ to an $A_n$ map lies in
\begin{equation}
[\Sigma^{n-2} A^{(n)},B] \cong B^{2-n}(A^{(n)}).
\end{equation}
Unitality implies that we can replace $A$ with $\bar{A}$.

We end the section by looking at how to endow a spectrum $M$ with an $A$-bimodule structure, given an $A_\infty$ structure on $A$. Given an $(A_{n-1},A)$-bimodule structure structure on $M$, we have to find maps
\begin{equation}
\xi_{n,i} : (K_n)_+ \sma A^{(i-1)} \sma M \sma A^{(n-i)} \lar M
\end{equation}
for $1 \leq i \leq n$. Each $\xi_{n,i}$ is determined by the $(A_{n-1},A)$-bimodule structure on $\partial K_n$, so the obstruction to defining $\xi_{n,i}$ lies in $M^{3-n}(A^{(i-1)} \sma M \sma A^{(n-i)})$.

As usual, unitality implies that the obstruction really lies in $M^{3-n}(\bar{A}^{(i-1)} \sma M \sma \bar{A}^{(n-i)})$. Similarly, the set of $(A_n,A)$-bimodule structures extending a given $(A_{n-1},A)$-bimodule structure is, if it is nonempty, a torsor over $\prod_{1 \leq i \leq n} M^{2-n}(\bar{A}^{(i-1)} \sma M \sma \bar{A}^{(n-i)})$. We see that there are generally more $A$-bimodule structures on $A$ than there are $A_\infty$ structures on $A$, because each of the maps $\xi_{n,i}$ for $1 \leq i \leq n$ can be chosen independently.

It is sometimes convenient to introduce another condition. There are two natural maps $A \sar M$, given by $A \cong A \sma_R R \sar A \sma_R M \sar M$ and $A \cong R \sma_R A \sar M \sma_R A \sar M$.  We can ask for these two maps to agree, and then for each composite map $(K_n)_+ \sma A^{(i-1)} \sma R \sma A^{(n-i)} \sar (K_n)_+ \sma A^{(i-1)} \sma M \sma A^{(n-i)} \sar M$ to be given by the $A_{n-1}$ structure on $A$ followed by the map $A \sar M$. In that case we can further reduce the obstruction to an obstruction in $M^{3-n}(\bar{A}^{(i-1)} \sma \bar{M} \sma \bar{A}^{(n-i)})$. We will call such an $A$-bimodule \emph{unital}.

For example, if $A=M=R/x$ with $R$ even commutative, we see that the set of unital $(A_n,A)$-bimodule structures on $A$ is a torsor over $(\pi_{n(d+2)-2} A)^n$. For convenience we will always assume that $M$ satisfies this condition, to make the theory of $A_\infty$ bimodules as closely related to the theory of $A_\infty$ structures as possible.

\section{Cyclic and cocyclic $A_n$ structures} \label{sec:cyclic}
It turns out that in certain cases, such as when $A=R/I$, the $A_n$ structure on $A$ controls certain hidden extensions of height $n-1$ in the canonical spectral sequences (\ref{homologySS} and \ref{cohomologySS}) calculating $\pi_*THH^R(A)$ and $\pi_*THH_R(A)$. For $THH^R(A)$ the extensions are trivial if and only if the maps $(K_n)_+ \sma A^{(n)} \sar A$ obtained by first cyclically permuting the $A$-factors and then using the $A_n$ structure are homotopic, in a sense we will make precise below. If these maps are homotopic we call the $A_n$ structure cyclic. For $THH_R(A)$ the cyclic permutations are replaced by the maps $A^{(n)} \sar A^{(n)}$ given by $(a_1,a_2,\ldots,a_n) \mapsto (a_2,\ldots,a_i,a_1,a_{i+1},\ldots,a_n)$, and if these maps are homotopic we call the $A_n$ structure cocyclic. We will come back to the $THH$-calculations in the next section. In this section we concentrate on developing the theory of cyclic and cocyclic $A_n$ structures.

Given an operad $P$ in a symmetric monoidal category $\C$ and a functor $R : {}^0 \! \bD C_P^{op} \sar \C$, recall from \cite[Definition 5.1]{An_cyclic} the definition of an $R$-trace on a pair $(A,M)$ consisting of a $P$-algebra $A$ and an $A$-bimodule $M$ in some symmetric monoidal $\C$-category $\D$ with target $B$ as a natural transformation $R \sar \E_{A,M,B}$ of functors. Here $\E_{A,M,B}(\{0,1,\ldots,n\})=Hom(M \otimes A^{\otimes n},B)$.

Also recall from \cite[Definition 5.3]{An_cyclic} the definition of an $R$-cotrace on a pair $(A,M)$ with source $B$ as a natural transformation $R \sar \tilde{\E}_{B,A,M}$, where $\tilde{\E}_{B,A,M}(\{0,1,\ldots,n\})=Hom(B \otimes A^{\otimes n},M)$. In particular, these definitions make sense if $P=\K$ is the associahedra operad and $R=\mathcal{W}$ is the cyclohedra, or if $P=\K_n$ and $R=\mathcal{W}_n$.

\begin{definition}
Let $A$ be an $A_n$-algebra, $1 \leq n \leq \infty$. We say that the $A_n$ structure on $A$ is cyclic if the identity map $A \sar A$ can be extended to a $\mathcal{W}_n$-trace.
\end{definition}

This means, roughly speaking, that the maps
\begin{equation}
(a_1,\ldots,a_n) \mapsto a_i \cdots a_n a_1 \cdots a_{i-1}
\end{equation}
for $1 \leq i \leq n$ are homotopic in a sufficiently nice way. To be more precise, a $\mathcal{W}_n$-trace on $A$ is a collection of maps $(W_m)_+ \sma A^{(m)} \sar A$ for $m \leq n$. The map $(W_n)_+ \sma A^{(n)} \sar A$ is determined on each of the $n$ copies of $K_n$ on the boundary of $W_n$ by the $A_n$ structure precomposed with a cyclic permutation, and given a $\mathcal{W}_{n-1}$-trace the map $(W_n)_+ \sma A^{(n)} \sar A$ is determined on all of $\partial W_n$. Thus a $\mathcal{W}_n$-trace on $A$ is a coherent choice of homotopies between the $n$ maps $(K_n)_+ \sma A^{(n)} \sar A$ obtained by first cyclically permuting the $A$-factors and then using the $A_n$ structure, plus some extra homotopies to glue these maps together.

\begin{figure}[h]
\epsfig{file=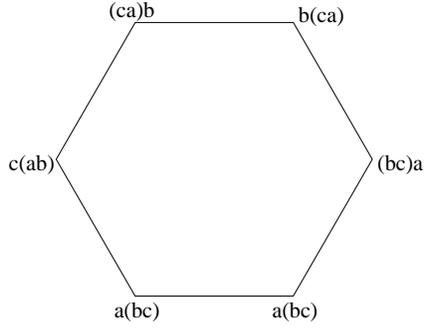,scale=0.7} \caption{The cyclohedron that has to be filled for a cyclic $A_3$ structure.} \label{fig:cyclicA3}
\end{figure}

A cyclic $A_2$ structure is the same as a homotopy commutative (and unital, but not necessarily homotopy associative) multiplication. A cyclic $A_3$ structure is an $A_3$ structure which is homotopy commutative, and such that for some choice of homotopy between the multiplication and its opposite, the natural map $(\partial W_3)_+ \sma A^{(3)} \sar A$ can be extended to all of $W_3$, as in Figure \ref{fig:cyclicA3}.

We can also make this definition for an $A$-bimodule $M$.

\begin{definition}
Let $A$ be an $A_n$-algebra and let $M$ be an $(A_n,A)$-bimodule. We say that the bimodule structure on $M$ is cyclic if the identity map $M \sar M$ can be extended to a $\mathcal{W}_n$-trace.
\end{definition}

Thus $M$ is a cyclic bimodule over the $A_n$-algebra $A$ if the maps
\begin{equation}
(m,a_2,\ldots, \\ a_n) \mapsto a_i \cdots a_n m a_2 \cdots a_{i-1}
\end{equation}
for $1 \leq i \leq n$ are homotopic in the same sense as above.

We make a similar definition for a cotrace:

\begin{definition}
Let $A$ be an $A_n$-algebra and let $M$ be an $(A_n,A)$-bimodule. We say that the $A_n$ structure on $A$ is cocyclic if the identity map $A \sar A$ can be extended to a $\mathcal{W}_n$-cotrace, and we say that the bimodule structure on $M$ is cocyclic if the identity map $M \sar M$ can be extended to a $\mathcal{W}_n$-cotrace.
\end{definition}

\begin{figure}[h]
\epsfig{file=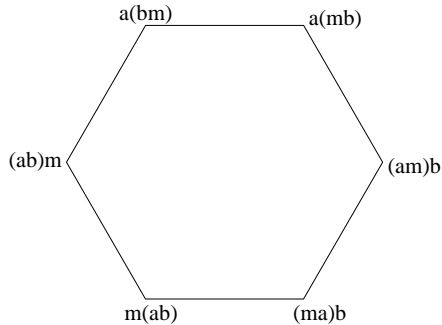,scale=0.7} \caption{The cyclohedron that has to be filled for a cocyclic $(A_3,A)$-bimodule structure.} \label{fig:cocyclicA3}
\end{figure}

Thus $M$ is cocyclic if the maps
\begin{equation}
(m,a_2,\ldots,a_n) \mapsto a_2 \cdots a_i m a_{i+1} \cdots a_n
\end{equation}
are homotopic. Note that the cyclic permutations from before have been replaced by a linear ordering of the $A$-factors while $M$ is allowed in any position. When $n=2$ cyclic and cocyclic mean the same thing, namely that the left and right action of $A$ on $M$ are homotopic, or that $M$ is homotopy symmetric. When $n=3$, $M$ is a cocyclic $(A_3,A)$-bimodule if the hexagon in Figure \ref{fig:cocyclicA3} can be filled. The boundary of this hexagon has the same shape as one of the two diagrams relating associativity to the twist map in the definition of a braided monoidal category.

Because a $\mathcal{W}$-trace is corepresented by the cyclic bar construction (\cite[Observation 5.2]{An_cyclic}), which in the category of spectra is topological Hochschild homology, and a $\mathcal{W}$-cotrace is represented by the cyclic cobar construction (\cite[Observation 5.5]{An_cyclic}), or topological Hochschild cohomology, this helps us study maps out of topological Hochschild homology and into topological Hochschild cohomology.

\begin{observation}
Let $A$ be an $A_\infty$ $R$-algebra and let $M$ be an $A$-bimodule. Then the natural map $M \sar sk_{n-1}THH^R(A;M)$ splits if and only if $M$, considered as an $(A_n,A)$-bimodule, is cyclic.

Similarly, the natural map $Tot^{n-1}THH_R(A;M) \sar M$ splits if and only if $M$, considered as an $(A_n,A)$-bimodule, is cocyclic.
\end{observation}

The notions of a trace and a cotrace are dual, in the following sense. If $M$ is an $A$-bimodule, the dual $DM=F_R(M,R)$ is again an $A$-bimodule in a natural way, and we can say the following:

\begin{proposition}
If $M$ is a cyclic $(A_n,A)$-bimodule, then $DM$ is a cocyclic $(A_n,A)$-bimodule. If $M$ is a cocyclic $(A_n,A)$-bimodule, then $DM$ is a cyclic $(A_n,A)$-bimodule.

If $M$ is dualizable in the sense that $DDM \cong M$, then $M$ is cyclic (cocyclic) if and only if $DM$ is cocyclic (cyclic).
\end{proposition}

\begin{proof}
This is clear, because the definition of a cotrace is dual to the definition of a trace. For example, for $W_3$ the homotopies $(ab)m \sim a(bm) \sim a(mb) \sim (am)b \sim (ma)b \sim m(ab) \sim (ab)m$ for $a,b \in A$ and $m \in M$ are dual to homotopies $f(ab) \sim (fa)b \sim b(fa) \sim (bf)a \sim a(bf) \sim (ab)f \sim f(ab)$ for $a,b \in A$ and $f \in DM$, and so a filling of the first hexagon is dual to a filling of the second hexagon.

The other case, and the if and only if statements when $M$ is dualizable, are similar.
\end{proof}

Because a $\mathcal{W}_n$-trace determines the map $(W_n)_+ \sma M \sma A^{(n-1)} \sar B$ on $\partial W_n$, we can play a similar game as with the $A_\infty$ obstruction theory. As before, let $\bar{A}$ be the cofiber of $R \sar A$, and let $\E_1^{s,t}=B^{-t}(M \sma \bar{A}^{(s)})$. Then, given a $\mathcal{W}_{n-1}$-trace, the obstruction to the existence of a $\mathcal{W}_n$-trace lies in $\E_1^{n-1,n-2}$. If the set of $\mathcal{W}_n$-traces is nonempty, it is a torsor over $\E_1^{n-1,n-1}$.

Similarly, let $\tilde{\E}_1^{s,t}=M^{-t}(B \sma \bar{A}^{(s)})$. Then, given a $\mathcal{W}_{n-1}$-cotrace, the obstruction to the existence of a $\mathcal{W}_n$-cotrace lies in $\tilde{\E}_1^{n-1,n-2}$.

The unreduced versions of $\E_1^{*,*}$ and $\tilde{\E}_1^{*,*}$ are graded cosimplicial abelian groups, and the following theorem follows in a similar way as Theorem \ref{AinftyobE_r}.

\begin{theorem}
Suppose we have a $\mathcal{W}_{n-1}$-trace extending a map $M \sar B$. Then, for $1 \leq r \leq n-1$, the obstruction to the existence of a $\mathcal{W}_n$-trace, while allowing the $\mathcal{W}_{n-r+1}$ through $\mathcal{W}_{n-1}$-trace to vary but fixing the $\mathcal{W}_{n-r}$-trace lies in $\E_r^{n-1,n-2}$.

Similarly, suppose we have a $\mathcal{W}_{n-1}$-cotrace extending a map $B \sar M$. Then, for $1 \leq r \leq n-1$, the obstruction to the existence of a $\mathcal{W}_n$-cotrace, while allowing the $\mathcal{W}_{n-r+1}$ through $\mathcal{W}_{n-1}$-cotrace to vary but fixing the $\mathcal{W}_{n-r}$-cotrace lies in $\tilde{\E}_r^{n-1,n-2}$.
\end{theorem}

For $A=R/I$ with $I=(x_1,\ldots,x_m)$ a regular ideal we will be especially interested in cotraces extending the natural map $R/x_i \sar R/I$. An important property of $R/I$ is that it is self dual (over $R$) up to a suspension. To be precise,
\begin{equation}
D(R/I) \simeq \Sigma^{-\Sigma(d_i+1)}R/I.
\end{equation}
The equivalence is not canonical, but a choice of regular sequence $(x_1,\ldots,x_m)$ generating $I$ dictates a choice of equivalence. This is clear from considering $R/x$ and $D(R/x)$, which are both $2$-cell $R$-modules with attaching map $x$.

Also note that the Bockstein $R/x \sar \Sigma^{d+1} R$ is dual to a map $\Sigma^{-d-1} R \sar D(R/x) \simeq \Sigma^{-d-1} R/x$, and this is precisely the unit map (desuspended) for $R/x$.

Thus the canonical map $R/x_i \sar R/I$ is dual to the map
\begin{equation}
R/I \sar \Sigma^{\Sigma_{j \neq i} (d_j+1)} R/x_i
\end{equation}
given by the Bockstein $\beta_j$ on each factor $R/x_j$ for $j \neq i$. We then have the following.

\begin{theorem}
There is a $\mathcal{W}_n$-cotrace extending the canonical map $R/x_i \sar R/I$ if and only if there is a $\mathcal{W}_n$-trace extending the canonical map $R/I \sar \Sigma^{\Sigma_{j \neq i} (d_j+1)} R/x_i$. In particular, an $A_n$ structure on $A=R/x_i$ is cyclic if and only if it is cocyclic.
\end{theorem}

\begin{proof}
For this we will need the spectral sequences (\ref{eq:THCSS}) and (\ref{eq:THHSS}), and the action described in Proposition \ref{prop:actionofSSs}. We have a $\mathcal{W}_n$-cotrace on $R/x_i \sar R/I$ if and only if $x_i$ acts trivially on $\pi_* THH_R(R/I)$ up to filtration $n-1$. Consider the pairing
\begin{equation}
THH_R(R/I) \sma_R THH^R(R/I) \sar THH^R(R/I),
\end{equation}
which is $R$-linear. Let $\gamma_I(\bar{q})$ denote some element $\gamma_{i_1}(\bar{q}_1) \cdots \gamma_{i_m}(\bar{q}_m)$ with $i_1+\cdots i_m=n-1$. If $x_i \gamma_I(\bar{q})=v \neq 0$, then $x_i(1,\gamma_I(\bar{q})) \mapsto v$, so $(x_i,\gamma_I(\bar{q})) \mapsto v$. But $x_i$ is in filtration $\geq n$ and $\gamma_I(\bar{q})$ is in filtration $n-1$, so in fact $(x_i,\gamma_I(\bar{q})) \mapsto 0$.
\end{proof}

It will be convenient to write down exactly what the cotrace obstructions look like in the cases we care about. Let $A=M=R/I$ and consider $\mathcal{W}_n$-cotraces from $B=R/x_i$ to $(A,M)$. Then
\begin{equation}
\tilde{\E}_1^{s,-*}=A^*(R/x_i \sma \bar{A}^{(s)}) \cong \Lambda_{A_*}(\alpha_i) \otimes_{A_*} \bar{E}_1^{s,-*},
\end{equation}
where $\bar{E}_1$ is the (reduced) $E_1$-term of the spectral sequence converging to $\pi_* THH_R(A)$. The $d_1$-differential on $\tilde{\E}_1$ is tensored up from the $d_1$-differential on $\bar{E}_1$, so
\begin{equation}
\tilde{\E}_2^{*,-*} \cong \Lambda_{A_*}(\alpha_i) \otimes_{A_*} A_*[q_1,\ldots,q_m].
\end{equation}
The obstruction to a $\mathcal{W}_n$-cotrace lies in $\tilde{\E}_2^{n-1,n-2}$, which has odd total degree. The only generator of $\tilde{\E}_2^{n-1,n-2}$ which has odd degree is $\alpha_i$, so the obstruction has to look like $\alpha_i r_i(q_1,\ldots,q_m)$ for some polynomial $r_i$ of degree $n-1$ in the $q_i$'s and total degree $d_i$. For notational convenience, we will write this obstruction as $r_i(q_1,\ldots,q_m)dq_i$. If $B=R/x_1 \vee \ldots \vee R/x_m$ then the obstruction looks like
\begin{equation}
\sum r_i(q_1,\ldots,q_m)dq_i
\end{equation}
for polynomials $r_1,\ldots,r_m$. We will see in Corollary \ref{cor:changebydf} that if we change the $A_n$ structure on $A=R/I$ by $f(q_1,\ldots,q_m)$ then this obstruction changes by $df=\sum \frac{\partial f}{\partial q_i} dq_i$.

\subsection*{The case $A=R/x$}
Now let us study the case $A=M=B=R/x$ with $R$ even and $x$ a nonzero divisor. We need to go into more detail about how to calculate the actual obstructions. We first do the $n=2$ case, following Strickland \cite{Str}. In this case we want to know how to determine whether or not an $A_2$ structure on $A$ is (co)cyclic, or in other words whether or not $A$ has a homotopy commutative multiplication.

Given a multiplication ($A_2$ structure) $\phi$ on $A=R/x$, we note that $\phi$ and $\phi^{op}$ agree on the bottom $3$ cells of $R/x \sma_R R/x$ regarded as a $4$-cell $R$-module, and following Strickland we define $c(\phi)$ by the equation
\begin{equation}
\phi^{op}-\phi=c(\phi) \circ (\beta \sma \beta).
\end{equation}
Clearly $\phi$ is homotopy commutative, or (co)cyclic, if and only if $c(\phi)=0$.

If $\phi'=\phi+v \circ (\beta \sma \beta)$ for some $v \in \pi_{2d+2} R/x$ then $(\phi')^{op}=\phi^{op}-v \circ (\beta \sma \beta)$, so $c(\phi')=c(\phi)-2v$. In particular, if $2$ is invertible in $\pi_* R/x$ then $R/x$ always has a cyclic $A_2$ structure. Strickland defines $\bar{c}(x)$ as the image of $c(\phi)$ in $\pi_{2d+2}R/(2,x)$.

Now the following proposition, except the last part about power operations, is clear:
\begin{proposition} \label{prop:Strickland} (Strickland, \cite[Proposition 3.1]{Str})
\begin{enumerate}
\item All products ($A_2$ structures) are associative ($A_3$), and unital.
\item The set of products on $R/x$ has a free transitive action of the group $R_{2d+2}/x$.
\item There is a naturally defined element $\bar{c}(x) \in \pi_{2d+2} R/(2,x)$ such that $R/x$ admits a commutative product if and only if $\bar{c}(x)=0$.
\item If so, the set of commutative products has a free transitive action of $ann(2,R_{2d+2}/x)=\{y \in R_{2d+2}/x \, | \, 2y = 0 \}$.
\item If $d \geq 0$ there is a power operation $\tilde{P} : R_d \sar R_{2d+2}/2$ such that $\bar{c}(x)=\tilde{P}(x) \,\, (mod \,\, 2,x)$ for all $x$.
\end{enumerate}
\end{proposition}

The power operation is constructed as follows. Consider the map $x^2 : \Sigma^{2d} R \sar R$. Because $R$ is $E_\infty$, we can extend this map over $\R P^{\infty}_{2d}$. By restricting to $\R P^{2d+2}_{2d} \simeq \Sigma^{2d} \R P^2$ we get a map $\Sigma^{2d} \R P^2_+ \sma R \sar R$, or an element $P(x)$ in $R^{-2d}(\R P^2_+ \sma R) \cong R_{2d} \oplus R_{2d+2}/2$. $\tilde{P}(x)$ is defined as the projection of $P(x)$ onto $R_{2d+2}/2$. The above proposition is also true when $d<0$, though some details in the proof would have to be changed.

The obstruction theory for cyclic $A_n$ structures for $n>2$ is somewhat harder. As before, given an $A_n$ structure $\phi=(\phi_2,\ldots,\phi_n)$ on $R/x$ where the $A_{n-1}$ structure is cyclic, we can define an obstruction $c_n(\phi_n) \in \pi_{n(d+2)-2} R/x$ to the $A_n$ structure being cyclic by factoring $\partial W_n \sma R/x^{(n)} \sar R/x$ as $\partial W_n \sma R/x^{(n)} \cong \Sigma^{n-2} R/x^{(n)} \overset{\beta \sma \ldots \sma \beta}{\lar} \Sigma^{n(d+2)-2}R \overset{c(\phi_n)}{\sar} R/x$.

Recall that the cyclohedron $W_n$ has $n$ copies of $K_n$ on its boundary. Now, if $\phi_n'=\phi_n+v \circ (\beta \sma \ldots \sma \beta)$, then the obstruction to the $A_n$ structure being cyclic changes on each of the $K_n$'s, and the net effect is that $c_n(\phi_n')=c_n(\phi_n)+nv$. We define $\bar{c}_n(x) \in \pi_{n(d+2)-2} R/(n,x)$ as the image of $c_n(\phi_n)$ for some $A_n$ structure extending a cyclic $A_{n-1}$ structure. It follows that if $n$ is invertible in $R_*/x$ then we can always find a cyclic $A_n$ structure extending a cyclic $A_{n-1}$ structure. If $n=p$ is not invertible, we make the following conjecture.

\begin{conjecture} \label{conj:powerop}
Suppose $n$ is invertible in $R_*$ for $n<p$ and $2(p-1) \, | \, d$. Then there is a power operation $\tilde{P} : R_d \sar R_{p(d+2)-2} R/p$ such that $\bar{c}_n(x)=\tilde{P}(x) \, \, (mod \,\, p,x)$ for all $x$.
\end{conjecture}

There certainly is such a power operation, constructed by extending $x^p : \Sigma^{pd} R \sar R$ over a skeleton of $B\Sigma_p$, the problem is identifying the power operation with the obstruction.

\subsection*{The Morava $K$-theories}
Conjecture \ref{conj:powerop} would generalize Proposition \ref{prop:Strickland}, and in particular this would show that there is no $\mathcal{W}_p$-trace (or cotrace) on $K(n)$. Instead we show this, and more, by finding $A_*$ comodule extensions in the B\"okstedt spectral sequence converging to $H_*(THH^S(k(n));\F_p)$. This will also supply us with the obstructions to $\mathcal{W}$-cotraces extending the natural maps $\hEn/v_i \lar K(n)$ for $0 \leq i \leq n-1$ ($v_0=p$). We will give the argument for odd primes; the $p=2$ case is similar after making the usual changes in the notation. Consider the connective Morava $K$-theory spectrum $k(n)$. From \cite{BaMa72} we know that
\begin{equation}
H_*(k(n);\F_p)=P(\bar{\xi}_i \, | \, i \geq 1) \otimes E(\bar{\tau}_i \, | \, i \neq n).
\end{equation}
The calculation of the $E^2$-term of the B\"okstedt spectral sequence converging to $H_*(THH^S(k(n));\F_p)$ is similar to the calculations found for example in \cite[\S 5]{AnRo}, and we get
\begin{equation}
E^2_{**}=H_*(k(n);\F_p) \otimes E(\sigma \bar{\xi}_i \, | \, i \geq 1) \otimes \Gamma(\sigma \bar{\tau}_i \, | \, i \neq n).
\end{equation}
Because there is no multiplication on $THH^S(k(n))$ this has to be interpreted additively only.

This $E^2$-term injects into the $E^2$-term for the corresponding spectral sequence for $H\Z/p$, so the differentials are induced by the corresponding differentials for $H\Z/p$. Thus there is a differential $d^{p-1}(\gamma_p(\sigma \bar{\tau}_i))=\sigma \bar{\xi}_{i+1}$, and the $E^p$-term looks like
\begin{equation}
E^p_{**}=H_*(k(n);\F_p) \otimes E(\sigma \bar{\xi}_{n+1}) \otimes P_p(\sigma \bar{\tau}_i \, | \, i \neq n).
\end{equation}
At this point the map of spectral sequences is no longer injective, so we cannot use this argument to say that the spectral sequence collapses. But we can say that there are no more differentials in low degrees:

\begin{proposition} \label{low_collapse}
The B\"okstedt spectral sequence converging to $H_*(k(n);\F_p)$ has no $d^n$ differentials for $n \geq p$ in degree less than $2p^{n+1}-1$.
\end{proposition}

\begin{proof}
By comparing with the B\"okstedt spectral sequence for $H\Z/p$, any differential has to hit something something which is in the kernel of the map $E^p_{**}(k(n)) \lar E^p_{**}(H\Z/p)$. The first element in the kernel is $\sigma \bar{\xi}_{n+1}$, which has degree $2p^{n+1}-1$.
\end{proof}

In particular, this shows that $\bigotimes_{0 \leq i <n} P_p(\sigma \bar{\tau}_i)$ survives to $E^\infty$.

\begin{remark}
If $n=1$, then one can show (\cite{AuRo2}) that the spectral sequence does collapse, by using that the map $\ell \lar k(1)$ makes the B\"okstedt spectral sequence for $k(1)$ into a module spectral sequence over the B\"okstedt spectral sequence for $\ell$, and using that $E^p_{**}(k(1))$ is generated as a module over $E^p_{**}(\ell)$ by classes in filtration $0 \leq i \leq p-1$. One could imagine a similar argument with $E^*_{**}(k(n))$ as a module over $E^*_{**}(BP\langle n \rangle)$ if $BP \langle n \rangle$ is at least $E_2$, though in this case the module generators are in filtration $0 \leq i \leq n(p-1)$.
\end{remark}

Recall that in the corresponding B\"okstedt spectral sequence for $H\Z/p$ there are multiplicative extensions $(\sigma \bar{\tau}_i)^p=\sigma \bar{\tau}_{i+1}$. Thus we find that 
\begin{equation}
\sigma ( \bar{\tau}_{n-1} (\sigma \bar{\tau}_{n-1})^{p-1})=\sigma \bar{\tau}_n 
\end{equation}
in $H_*(THH^S(H \Z/p),\F_p)$, and more generally
\begin{equation}
\sigma( \bar{\tau}_i (\sigma \bar{\tau}_i)^{p-1} \cdots (\sigma \bar{\tau}_{n-1})^{p-1})=\sigma \bar{\tau}_n.
\end{equation}
We use this to prove that there are $A_*$ comodule extensions in the B\"okstedt spectral sequence for $k(n)$.

\begin{proposition} \label{Acomodkn}
Let $x_i=(\sigma \bar{\tau}_i)^{p-1} \cdots (\sigma \bar{\tau}_{n-1})^{p-1}$ in $H_*(THH^S(k(n));\F_p)$. Then the $A_*$ comodule action on $\bar{\tau}_i x_i$ is given by
\begin{equation}
\nu(\bar{\tau}_i x_i)=1 \otimes \bar{\tau}_i x_i+\sum \bar{\tau}_j \otimes \bar{\xi}_{i-j}^{p^j}x_i-\sum \bar{\tau}_j \otimes \bar{\xi}_{n-j}^{p^j}.
\end{equation}

All the classes in $\bigotimes_{0 \leq i < n} P_p(\sigma \bar{\tau}_i)$ are $A_*$ comodule primitive, and together with the natural $A_*$ comodule structure on $H_*(k(n))$ this determines the $A_*$ comodule structure on $H_*(THH(k(n));\F_p)$ up to degree $2p^{n+1}-1$.
\end{proposition}

\begin{proof}
Consider the commutative diagram
\begin{equation}
\xymatrix{ H_*(THH^S(k(n));\F_p) \ar[r]^-\sigma \ar[d] & H_{*+1}(THH^S(k(n));\F_p) \ar[d] \\
H_*(THH^S(H\Z/p);\F_p) \ar[r]^-\sigma & H_{*+1}(THH^S(H\Z/p);\F_p) }
\end{equation}
The classes in question all survive to $E^\infty_{**}$ by Proposition \ref{low_collapse}, and because $\sigma(\bar{\tau}_i x_i)=0$ in $H_*(THH^S(k(n));\F_p)$,  we conclude that the image of $\bar{\tau}_i x_i$ in $H_*(THH^S(H\Z/p);\F_p)$ has to be in the kernel of $\sigma$. But $\sigma(\bar{\tau}_i x_i)=\sigma \bar{\tau}_n$ in $H_*(THH^S(H\Z/p);\F_p)$, and the image is given by the element with the same name in the B\"okstedt spectral sequence for $H\Z/p$ modulo lower filtration. Thus $\bar{\tau}_i x_i$ in $H_*(THH^S(k(n));\F_p)$ has to map to $\bar{\tau}_i x_i$ minus something in lower filtration which also maps to $\sigma \bar{\tau}_n$ under $\sigma$. The only elements in lower filtration that map to $\sigma \bar{\tau}_n$ are $\bar{\tau}_n$ and $\bar{\tau}_j x_j$ for $j>i$. If necessary we can adjust $\bar{\tau}_i x_i$ by adding elements in lower filtration in the B\"okstedt spectral sequence for $k(n)$ so $\bar{\tau}_i x_i$ maps to $\bar{\tau}_i x_i-\bar{\tau}_n$.

Because the map from $H_*(THH^S(k(n));\F_p)$ to $H_*(THH^S(H\Z/p);\F_p)$ is a map of $A_*$ comodules, it follows that the $A_*$ comodule action on $\bar{\tau}_i x_i$ is as claimed. The claim about all the classes in $\bigotimes_{0 \leq i < n} P_p(\sigma \bar{\tau}_i)$ being primitive follows immediately by using that
\begin{equation}
H_*(THH^S(k(n));\F_p) \lar H_*(THH^S(H\Z/p);\F_p)
\end{equation}
is injective in low degrees.
\end{proof}

Now let us see what happens in the Adams spectral sequence with $E_2$-term $Ext_{A_*}(\F_p,H_*(THH^S(k(n));\F_p))$ converging to $\pi_*THH^S(k(n))$. (Here $Ext$ means $Ext$ of comodules, as opposed to in the topological Hochschild cohomology spectral sequence.) Because everything is concentrated in Adams filtration $0$ and $1$ in low degrees, we can run the whole Adams spectral sequence through a range of dimensions.

\begin{theorem} \label{K(n)obstructions}
The classes $\bigotimes_{0 \leq i < n}P_p(\sigma \bar{\tau}_i)$ and $v_n$ all give rise to corresponding nontrivial classes in $\pi_*THH^S(k(n))$. Moreover, there are relations
\begin{equation}
v_i (\sigma \bar{\tau}_i)^{p-1} \cdots (\sigma \bar{\tau}_{n-1})^{p-1}=v_n.
\end{equation}
\end{theorem}

\begin{proof}
First of all, since all the classes in $\bigotimes_{0 \leq i < n}P_p(\sigma \bar{\tau}_i)$ are primitive, we get corresponding classes in filtration $0$ in the Adams spectral sequence. Also, because $\bar{\tau}_n$ is missing from $H_*(k(n);\F_p)$ we get a class $v_n$ in filtration $1$. There are no classes in higher filtration in these degrees, so the classes $\bigotimes_{0 \leq i <n}P_p(\sigma \bar{\tau}_i)$ and $v_n$ all survive to $\pi_* THH^S(k(n))$.

Recall, e.g.~from \cite[p.~63]{Ra04} that $v_n$ is represented by $-\sum \bar{\tau}_i \otimes \bar{\xi}_{n-i}^{p^i}$ in the cobar complex for $p$ odd, with a similar formula for $p=2$. This also implies that $v_i (\sigma \bar{\tau}_i)^{p-1} \cdots (\sigma \bar{\tau}_{n-1})^{p-1}$ is represented by $-\sum \bar{\tau}_j \otimes \bar{\xi}_{i-j}^{p^j} (\sigma \bar{\tau}_i)^{p-1} \cdots (\sigma \bar{\tau}_{n-1})^{p-1}$.

From the $A_*$ comodule structure we found in Proposition \ref{Acomodkn}, we find that the expressions representing $v_i (\sigma \bar{\tau}_i)^{p-1} \cdots (\sigma \bar{\tau}_{n-1})^{p-1}$ and $v_n$ are homologous, so the two expressions have to be equal in $\pi_*THH^S(k(n))$.
\end{proof}

\section{Calculations of $THH(A)$ for $A=R/I$}
In this section we attempt to calculate $THH_R(R/I)$ and $THH^R(R/I)$, where as usual $R$ is even commutative and $I=(x_1,\ldots,x_m)$ is a finitely generated regular ideal with $|x_i|=d_i$. Some of the results in this section also hold when $I$ is infinitely generated.

By Proposition \ref{prop:AAopring} we know that $\pi_* A \sma_R A^{op} \cong \Lambda_{A_*}(\alpha_1,\ldots,\alpha_m)$ with $|\alpha_i|=d_i+1$ and by the $Ext$ calculation in the proof of Corollary \ref{cor:evenquotientainfty} we find that the $E_2$-term of the spectral sequence converging to $\pi_* THH_R(A)$ looks like
\begin{equation}
E_2^{*,*}=A_*[q_1,\ldots,q_m] \Longrightarrow \pi_* THH_R(A), \label{eq:THCSS}
\end{equation}
with $q_i$ in bidegree $(1,-d_i-1)$, or total degree $-d_i-2$.

A similar calculation shows that the $E_2$-term of the corresponding spectral sequence converging to $\pi_* THH^R(A)$ looks like
\begin{equation}
E^2_{*,*}=\Gamma_{A_*}[\bar{q}_1,\ldots,\bar{q}_m] \Longrightarrow \pi_* THH^R(A), \label{eq:THHSS}
\end{equation}
where $\Gamma$ denotes a divided power algebra. In this case $\bar{q}_i$ is in bidegree $(1,d_i+1)$, or total degree $d_i+2$. The first spectral sequence is a spectral sequence of algebras, while the second spectral sequence has to be interpreted additively only unless the multiplication on $A$ is more structured than just $A_\infty$. For example, if $A$ is $E_n$ for some $n \geq 2$, then by \cite{BaMa} or \cite{FiVo} $THH^R(A)$ is $E_{n-1}$ and this also gives a multiplication on the spectral sequence.

\begin{proposition} \label{prop:actionofSSs}
The $E_2$-term of the spectral sequence (\ref{eq:THCSS}) acts on the $E^2$-term of the spectral sequence (\ref{eq:THHSS}) by $q_i \cdot \gamma_j(\bar{q}_i)=\gamma_{j-1}(\bar{q}_i)$ and $q_i \cdot \gamma_j(\bar{q}_k)=0$ for $k \neq i$.
\end{proposition}

Here $\gamma_0(\bar{q}_i)=1$ and $\gamma_j(\bar{q}_i)=0$ for $j<0$.

\begin{proof}
This follows because algebraically a divided power algebra is dual to a polynomial algebra, and $Tor$ and $Ext$ are dual in this case.
\end{proof}

Both spectral sequences collapse at the $E_2$-term, and we are now in a position to look for hidden extensions in the two spectral sequences. The extensions are all of the following form. Each $x_i \in \pi_{d_i} R$ acts trivially on $E_{\infty}$, but might act nontrivially on $\pi_* THH_R(A)$ and $\pi_* THH^R(A)$. In the first spectral sequence an element is well defined modulo higher filtration, so multiplication by $x_i$ increases the filtration, while in the second spectral sequence an element is well defined modulo lower filtration, so multiplication by $x_i$ decreases the filtration.

\subsection*{Height $1$ extensions} We will consider $A=R/x$ with $|x|=d$ first. In this case we have spectral sequences
\begin{eqnarray}
E_2^{*,*}=A_*[q] \Longrightarrow \pi_* THH_R(A), \label{eq:THCSSn=1} \\
E^2_{*,*}=\Gamma_{A_*}[\bar{q}] \Longrightarrow \pi_* THH^R(A). \label{eq:THHSSn=1}
\end{eqnarray}

Baker and Lazarev \cite{BaLa} have one result in this direction. They consider $R=KU$ and $A=KU/2$ and go on to calculate $THH_{KU}(KU/2)$. Their main tool is the following piece of Morita theory, a kind of double centralizer theorem which is an easy consequence of the theory developed in \cite{DwGr02}.
\begin{theorem} (Baker and Lazarev, \cite{BaLa}) \label{thm:doublecentralizer}
For a finite cell $R$-module $A$, the natural map
\begin{equation}
R \sar F_{F_R(A,A)}(A,A)
\end{equation}
is an $A$-localization. (Here $F(-,-)$ denotes the \emph{derived} function spectrum.)
\end{theorem}
In particular, $F_{F_{KU}(KU/2,KU/2)}(KU/2,KU/2) \simeq KU^\wedge_2$, and Baker and Lazarev proved that $THH_{KU}(KU/2) \simeq KU^\wedge_2$ by comparing $KU/2 \sma_{KU} KU/2^{op}$ with $F_{KU}(KU/2,KU/2)$.
\begin{theorem} (Baker and Lazarev, \cite[Theorem 3.1]{BaLa})
\begin{equation}
THH_{KU}(KU/2) \simeq KU^\wedge_2,
\end{equation}
where $KU^\wedge_2$ is the $2$-complete $K$-theory spectrum.
\end{theorem}

We will discuss the proof here because we want to compare it to our more general method. Consider the map $KU/2 \sma_{KU} KU/2^{op} \sar F_{KU}(KU/2,KU/2)$ adjoint to the action map of $KU/2 \sma_{KU} KU/2^{op}$ on $KU/2$. On homotopy this looks like
\begin{equation}
\Lambda_{KU/2_*}(\alpha) \sar \Lambda_{KU/2_*}(Q).
\end{equation}
Here $|\alpha|=1$ and $|Q|=-1$, so $\alpha \mapsto vQ$ for some $v \in \pi_2 KU/2$. The action map $KU/2 \sma_{KU} KU/2^{op} \sma_{KU} KU/2 \sar KU/2$ uses the twist map, so not surprisingly (\cite[Proposition 2.12]{BaLa}) $v$ measures the noncommutativity of $KU/2$, and in fact $v=c(\phi)=u \in \pi_2 KU/2$ is nontrivial. It follows that $KU/2 \sma_{KU} KU/2^{op} \sar F_{KU}(KU/2,KU/2)$ is a weak equivalence, and thus
\begin{equation}
F_{F_{KU}(KU/2,KU/2)}(KU/2,KU/2) \lar F_{KU/2 \sma_{KU} KU/2^{op}}(KU/2,KU/2)
\end{equation}
is a weak equivalence, proving the theorem.

It is clear that a generalization of their proof will tell us exactly when $THH_R(A)$ is the $A$-localization of $R$, but otherwise this method will not tell us much more about $THH_R(A)$. If we consider the spectral sequence (\ref{eq:THCSSn=1}), then what Baker and Lazarev really do is identifying an additive extension in the spectral sequence. While multiplication by $2$ acts trivially on $E_2=E_\infty$, in $\pi_* THH_{KU}(KU/2)$ multiplication by $2$ acts nontrivially. By abuse of notation, let $q$ also denote a representative of $q$ in $\pi_* THH_{KU}(KU/2)$. Then $2 \cdot 1 \equiv u q$ modulo filtration $2$ and higher, which is enough to conclude that $\pi_* THH_{KU}(KU/2) \cong \Z_p[u,u^{-1}]$ and that $KU \sar THH_{KU}(KU/2)$ is a $2$-completion.

Let us also consider $THH^{KU}(KU/2)$ before we return to the more general situation. In this case the spectral sequence looks almost identical, except the additive extension will decrease the filtration rather than increasing it. If we let $\gamma_i(\bar{q})$ also denote a representative of $\gamma_i(\bar{q})$ in $\pi_* THH^{KU}(KU/2)$, then $2 \cdot \gamma_i(\bar{q})=u \gamma_{i-1}(\bar{q})$ modulo filtration $i-2$. One way to see this is as follows:

\begin{lemma} \label{lem:THHTHCduality}
Suppose there is an extension $x=\sum a_i q^i$ in the spectral sequence (\ref{eq:THCSSn=1}). Then there are extensions
\begin{equation}
x \gamma_n(\bar{q})=\sum a_i \gamma_{n-i}(\bar{q})
\end{equation}
in the spectral sequence (\ref{eq:THHSSn=1}).
\end{lemma}

\begin{proof}
Consider the pairing $E_2^{**} \otimes_{R_*} E^2_{**} \sar E^2_{**}$ and the corresponding pairing $THH_R(A) \sma_R THH^R(A) \sar THH^R(A)$. Then
\begin{equation}
(x,\gamma_n(\bar{q}))=(\sum a_i q^i,\gamma_n(\bar{q})),
\end{equation}
which maps to $\sum a_i \gamma_{n-i}(\bar{q})$.
\end{proof}

Thus we get the following description of $THH^{KU}(KU/2)$:
\begin{theorem}
Topological Hochschild homology of $KU/2$ is given by
\begin{equation}
THH^{KU}(KU/2) \simeq KU[1/2]/2KU.
\end{equation}
\end{theorem}

\begin{proof}
By the lemma, $1 \in \pi_0 THH^{KU}(KU/2)$ is $2$-divisible, and the same argument shows that it is infinitely $2$-divisible. We can choose $1/2$ to be in filtration $1$, and from the spectral sequence it then follows that
\begin{equation}
\pi_n THH^{KU}(KU/2) \cong \Z/2^\infty
\end{equation}
for $n$ even, and the result follows.
\end{proof}

Note that while $THH_{KU}(KU/2)$ is automatically $KU/2$-local, the same is not true for $THH^{KU}(KU/2)$. If we localize (or $p$-complete) $THH^{KU}(KU/2)$, we end up with $\Sigma KU_2^\wedge$, indicating an interesting kind of duality between topological Hochschild homology and cohomology.

A straightforward generalization of Baker and Lazarev's method works equally well to determine when $THH_R(A) \simeq R_A$, the $A$-localization of $R$, for $A=R/I$.

For a multiplication ($A_2$ structure) $\phi$ on $A=R/I$ which is homotopy associative (can be extended to an $A_3$, and thus to an $A_\infty$ structure) we define an $n \times n$ matrix $C(\phi)$ as follows. If $\phi^0$ is given by smashing together multiplications $\phi^i$ on $R/x_i$ for each $i$ we set $c_{ii}(\phi^0)=c(\phi^i)$ and $c_{ij}(\phi^0)=0$ for $i \neq j$. If $\phi=\phi^0 \prod_{i \neq j}(1 \sma 1+v_{ij} Q_i \sma Q_j)$ we set $c_{ij}(\phi)=-v_{ij}-v_{ji}$. Thus $C(\phi)=0$ if and only if $\phi$ is homotopy commutative, so we can say that $C(\phi)$ measures the noncommutativity of the multiplication.

\begin{remark} \label{rem:autos}
If $v \in \pi_{d_i+d_j+2} A$, then we can construct an endomorphism $vQ_i Q_j$ of $A$, as follows:
\begin{equation}
A \arrow{Q_j} \Sigma^{-d_j-1} A \arrow{Q_i} \Sigma^{-d_i-d_j-2} A \arrow{v} A.
\end{equation}
Let $e=id+vQ_iQ_j$. Then $e$ is an automorphism of $A$ with inverse $e^{-1}=id-vQ_iQ_j$. If we have a multiplication $\phi$ on $A$, we can conjugate by $e$ to get a new multiplication $\phi^e$ defined by
\begin{equation}
A \sma A \arrow{e^{-1} \sma e^{-1}} A \sma A \arrow{\phi} A \arrow{e} A.
\end{equation}
One can check that $\phi^e=\phi(1 \sma 1-vQ_i \sma Q_j)(1 \sma 1+vQ_j \sma Q_i)$, so $C(\phi^e)=C(\phi)$. If $\phi=\phi^0 \prod_{i \neq j}(1 \sma 1+v_{ij} Q_i \sma Q_j)$ and $\phi'=\phi^0 \prod_{i \neq j} (1 \sma 1+v'_{ij} Q_i \sma Q_j)$ for the same $\phi_0$, then $C(\phi)=C(\phi')$ if and only if $\phi'$ is obtained from $\phi$ by conjugating by an endomorphism. Thus, at least off the diagonal the matrix $C(\phi)$ determines $\phi$ up to conjugation.
\end{remark}

\begin{theorem}
Suppose $A=R/I$ and let $\phi$ be an $A_2$ structure on $A$ which extends to an $A_\infty$ structure. If $C(\phi)$ is invertible we get
\begin{equation}
THH_R(A) \simeq R_A,
\end{equation}
the $A$-localization of $R$.
\end{theorem}

\begin{proof}
The map $A \sma_R A^{op} \sar F_R(A,A)$ sends $\alpha_i$ to $\sum c_{ij} \beta_j$, so we get an equivalence $A \sma_R A^{op} \overset{\simeq}{\sar} F_R(A,A)$ if and only if $C(\phi)$ is invertible. The result then follows from Theorem \ref{thm:doublecentralizer}.

Alternatively, we can identify the extensions in the spectral sequence (\ref{eq:THCSS}). Let $B=R/x_i$, and consider the obstruction to a $\mathcal{W}_2$-cotrace from $B$ to $(A,A)$. The obstruction is $\sum_j c_{ij} dx_i$. On the other hand, the obstruction is also the obstruction to the existence of the dotted arrow in the diagram
\begin{equation}
\xymatrix{R \ar[r]^{x_i} & R \ar[d] \ar[r] & R/x_i \ar@{.>}[ld] \\ & Tot^1 THH_R(A) & }
\end{equation}
which is exactly the nontriviality of $x_i \in \pi_{d_i} Tot^1 THH_R(A)$. Thus $x_i \equiv \sum_j c_{ij} q_j$ modulo filtration $2$ and higher. By the Weierstrass preparation theorem this is enough to conclude that when $C(\phi)$ is invertible then $\pi_* THH_R(A) \\ \cong (R_*)_I^\wedge$.
\end{proof}

Since this result holds whenever $C(\phi)$ is invertible, it holds generically in some sense. On the other hand, it can happen that each $c_{ij} \in \pi_{d_i+d_j+2} R/I$ is $0$ for degree reasons. If $R$ is $2$-periodic we have no degree considerations, and we can say that if $2$ is invertible in $A_*$ or if $n \geq 2$ then there exists a multiplication $\phi$ on $A$ with $C(\phi)$ invertible. In the characteristic $2$ case we use Proposition \ref{prop:Strickland}, part $5$ to determine whether or not $THH_R(R/x)$ is weakly equivalent to $R_{R/x}$.

The companion theorem for topological Hochschild homology is as follows.

\begin{theorem} \label{thm:Cinvertible}
Suppose $A=R/I$ and let $\phi$ be an $A_2$ structure on $A$ which extends to an $A_\infty$ structure. If $C(\phi)$ is invertible we get
\begin{equation}
THH^R(A) \simeq R[I^{-1}]/IR.
\end{equation}
\end{theorem}

\begin{proof}
As in Lemma \ref{lem:THHTHCduality}, $1 \in \pi_0 THH^R(A)$ is infinitely $x_i$-divisible for each $x_i$, with $1/x_i$ in filtration $1$.
\end{proof}

To proceed further it is useful to consider $THH(A;M)$ with $A=M=R/I$ as $R$-modules, but with a possibly different $A$-bimodule structure on $M$. We will assume all our $(A_n,A)$-bimodule structures are unital, in the sense described at the end of section \ref{sec:Ainftyobtheory}. Let $\xi_2=(\xi_{2,1},\xi_{2,2})$ be an $(A_2,A)$-bimodule structure on $M$. Here $\xi_{2,1} : M \sma A \sar M$ and $\xi_{2,2} : A \sma M \sar M$.

Now, it is possible to vary $\xi_{2,1}$ and $\xi_{2,2}$ independently. For example, we can get a new $(A_2,A)$-bimodule structure on $M$ by letting $\xi_{2,1}'=\xi_{2,1}(1 \sma 1+v_{ij} Q_i \sma Q_j)$ and not changing $\xi_{2,2}$. We can define a matrix $C(\xi_2)$ which measures how nonsymmetric the $(A_2,A)$-bimodule structure is in a similar way as before. If $A=M=R/x$, we define $c(\xi)$ by the equation
\begin{equation}
\xi_{2,1} \circ \tau - \xi_{2,2}=c(\xi)(\beta \sma \beta).
\end{equation}
If $\xi_2$ is obtained by smashing together $(A_2,R/x_i)$-bimodule structures $\xi_2^i$ on each $R/x_i$ we set $c_{ii}(\xi_2)=c(\xi_2^i)$. If $\xi_2'=(\xi_{2,1}',\xi_{2,2})$ where $\xi_{2,1}'$ is obtained from $\xi_{2,1}$ by $\xi_{2,1}'=\xi_{2,1}(1 \sma 1+v_{ij} Q_i \sma Q_j)$ we set $c_{ij}(\xi_2)=c_{ij}(\xi_2)-v_{ij}$. If $\xi_{2,2}$ is changed we adjust $c_{ji}$ in the same way.

It is clear that for any, not necessarily symmetric, $n \times n$ matrix $C$ with $c_{ij} \in \pi_{d_i+d_j+2} M$ there is an $(A_2,A)$-bimodule structure $\xi_2$ on $M$ with $C(\xi_2)=C$. It is also clear that $C(\xi_2)$ controls the height $1$ extensions in the $THH$ spectral sequences, so by choosing the $(A_2,A)$-bimodule structure appropriately we can get any extension which is possible for degree reasons.

\begin{theorem}
Suppose $A=M=R/I$, and fix an $A_\infty$ structure on $A$. Suppose there exists a matrix $C=\{c_{ij}\}$ with $c_{ij} \in \pi_{d_i+d_j+2} M$ which is invertible in $M_*$. Then there exists an $(A_\infty,A)$-bimodule structure on $M$ such that
\begin{equation}
THH_R(A;M) \simeq R_A
\end{equation}
and
\begin{equation}
THH^R(A;M) \simeq R[I^{-1}]/IR.
\end{equation}
\end{theorem}

\begin{proof}
Choose $\xi_2$ such that $C(\xi_2)=C$ is invertible.
\end{proof}

\subsection*{Height $n$ extensions}
Now we claim that this works equally well for $(A_n,A)$-bimodule structures on $M$ for any $n \geq 2$. The essential point is that if we perturb the $(A_n,A)$-bimodule structure while fixing the $(A_{n-1},A)$-bimodule structure, then all the height $n-2$ extensions remain the same and the height $n-1$ extensions change in a predictable way.

Suppose we have two $(A_n,A)$-bimodule structures $\xi=(\xi_2,\ldots,\xi_n)$ and $\xi'=(\xi_2',\ldots,\xi_n')$ on $M=A=R/I$ with $\xi_i=\xi'_i$ for $i<n$. Consider $\xi_n=\{\xi_{n,i} \, | \, 1 \leq i \leq n\}$ and $\xi'_n=\{\xi'_{n,i} \, | \, 1 \leq i \leq n\}$. Let $J=(j_1,\ldots,j_n)$ and write $Q_J$ for $Q_{j_1} \sma \ldots \sma Q_{j_n}$. Also write $q^J$ for $q_{j_1} \cdots q_{j_n}$. The obstruction theory implies that we can write
\begin{equation}
\xi'_{n,i}=\xi_{n,i} \prod_J \big(1 \sma \ldots \sma 1+v_{J,i} Q_J \big).
\end{equation}

\begin{theorem}
Let $A=M=R/I$ and suppose $M$ has two $(A_n,A)$-bimodule structures $\xi$ and $\xi'$ as above. Then the extension of $x_k$ to filtration $n-2$ is the same for the two $(A_n,A)$-bimodule structures, and the extension of $x_k$ to filtration $n-1$ differs by a sum over $i$ of all $J$ such that $j_i=k$ of $v_{J,i} q^J/q_k$.
\end{theorem}

\begin{proof}
This is mostly bookkeeping. Given $\xi$ and $\xi'$ as above, the difference between the degree $n-1$ extension for the two bimodule structures is given by the obstruction to the existence of a certain map $(W_n)_+ \sma R/x_k \sma A^{(n-1)} \sar M$, given by
\begin{equation}
(W_n)_+ \sma R/x_k \sma A^{(n-1)} \lar (W_n)_+ \sma M \sma A^{(n-1)} \lar M,
\end{equation}
where the first map is the canonical map $R/x_k \sar R/I$ and the second map is given by $\xi_{n,i}-\xi'_{n,i}$ on the $K_n$ on the boundary of $W_n$ corresponding to $(K_n)_+ \sma A^{(i-1)} \sma M \sma A^{(n-i-1)} \sar M$, and the trivial map on the rest of $\partial W_n$ because the $(A_{n-1},A)$-bimodule structures determined by $\xi$ and $\xi'$ agree. Using that $W_n \cong D^{n-1}$ and that the map on the boundary has to factor through $R/x \arrow{Q_k} \Sigma^{d_k+1} R$ the obstruction is given by a map
\begin{equation}
\Sigma^{n-2} R/x_k \sma A^{(n-1)} \arrow{Q_k \sma 1} \Sigma^{d_i+(n-1)}A^{(n-1)} \lar M,
\end{equation}
or an element in $M^{-d_i-(n-1)}(A^{(n-1)})$. On the $i$'th copy of $K_n$ the map $(K_n)_+ \sma A^{(i-1)} \sma R/x_k \sma A^{(n-i-1)} \sar M$ is given by a sum over all $J$ of
\begin{equation}
(K_n)_+ \sma A^{(i-1)} \sma R/x_k \sma A^{(n-i-1)} \arrow{\epsilon \sma 1} \Sigma^{d_k+1} (K_n)_+ \sma A^{(n-1)} \lar M,
\end{equation}
where $\epsilon=Q_k$ if $j_i=k$ and $0$ otherwise, and the second map is $v_{J,i} q^J/q_k$. This implies the result.
\end{proof}

It follows from this theorem that by adjusting the $(A_\infty,A)$-bimodule structure on $M$ we have a lot of choice for $\pi_* THH_R(A;M)$.

\begin{corollary} \label{cor:allTHHforM}
Given power series $f_i=\sum a_{iJ} q^J$ for $i=1,\ldots,m$ in $R_*[[q_1,\ldots,q_m]]$ with no constant term of total degree $d_1,\ldots,d_m$ respectively, there is an $(A_\infty,A)$-bimodule structure on $M$ such that
\begin{equation}
\pi_* THH_R(A;M) \cong R_*[[q_1,\ldots,q_m]]/(x_1-f_1,\ldots,x_m-f_m).
\end{equation}
For this $(A_\infty,A)$-bimodule structure we have
\begin{equation}
\pi_* THH^R(A;M) \cong \Gamma_{R_*}[\bar{q}_1,\ldots,\bar{q}_m]/(x_i \gamma_I(\bar{q})-\sum a_{iJ}\gamma_{I-J}(\bar{q})).
\end{equation}
Here $I$ and $J$ run over indices $I=(i_1,\ldots,i_m)$ and $J=(j_1,\ldots,j_m)$. Also, $q^J$ means $q_1^{j_1} \cdots q_m^{j_m}$, and $\gamma_I(\bar{q})$ means $\gamma_{i_1}(\bar{q}_1) \cdots \gamma_{i_m}(\bar{q}_m)$.
\end{corollary}

For generic power series $f_1,\ldots,f_m$ we see that $\pi_*THH_R(A;M)$ is a finite extension of $\pi_* R_A \cong (R_*)_I^\wedge$, and $\pi_*THH^R(A;M)$ is a finite free $R_*[I^{-1}]/IR_*$-module.  It also has the following consequence, now for $A=M$:

\begin{corollary} \label{cor:changebydf}
Suppose we change the $A_n$ structure on $A$ by some $f$ in $A_*[q_1,\ldots,q_m]$ of degree $n$ in the $q_i$'s. Then the extension of $x_k$ to filtration $n-1$ in $\pi_* THH_R(A)$ with these two $A_n$ structures (and any $A_\infty$ extensions of these) changes by $\frac{df}{dq_k}$.
\end{corollary}

This should be compared with Lazarev's Hochschild cohomology calculations in \cite{La03b}. He proved that if $R$ is an even commutative graded ring and $A$ is a $2$-cell DG $R$-module $A=\{\Sigma^{d+1} R \overset{x}{\sar} R\}$ with differential $x$ for some nonzero divisor $x$, then the moduli space of $A_\infty$ structures on $A$ can be identified with the set of power series $f(q)=xq+a_2 q^2+\ldots$ in $H_*(A)[[q]] \cong R/x[[q]]$. In this case, Lazarev proved (\cite[Proposition 7.1]{La03b}) that
\begin{equation}
HH_R^*(A) \cong H_*(A)[[q]]/(f'(q)).
\end{equation}
Thus the coefficient $a_n$ describing the $A_n$ structure contributes to an extension in the canonical spectral sequence, where multiplication by $x$ sends $1$ to $n a_n q^{n-1}$ in filtration $n-1$. But an $A_n$ structure is really $n$ identical maps $\xi_{n,i}$ making up an $(A_n,A)$-bimodule structure, and the coefficient $n$ comes from the sum of these $n$ maps.

\subsection*{The Morava $K$-theories}
There are two cases to consider, the $2(p^n-1)$-periodic Morava $K$-theory $K(n)=\hEn/(p,v_1,\ldots,v_{n-1})$ and the $2$-periodic Morava $K$-theory $K_n=E_n/(p,v_1,\ldots,v_{n-1})$. For $K(n)$ most of the (co)trace obstructions are zero for degree reasons, and we conjecture that $THH(K(n))$ is independent of the $A_\infty$ structure. On the other hand, $THH(K_n)$ varies considerably over the moduli space of $A_\infty$ structures.

We start by connecting the relation
\begin{equation}
v_i (\sigma \bar{\tau}_i)^{p-1} \cdots (\sigma \bar{\tau}_{n-1})^{p-1}=v_n
\end{equation}
in $\pi_* THH^S(k(n))$ from Theorem \ref{K(n)obstructions} to the $\mathcal{W}$-trace obstruction theory.

\begin{proposition} \label{prop:sigmatau=qbar}
Under the composite $THH^S(k(n)) \sar THH^S(K(n)) \sar THH^{\widehat{E(n)}}(K(n))$, the class $\sigma \bar{\tau}_i$ in $\pi_* THH^S(k(n))$ maps to the class $\bar{q}_i$ in $\pi_* THH^{\widehat{E(n)}}(K(n))$.
\end{proposition}

\begin{proof}
The key fact is that the exterior generator $\alpha_i$ in $\pi_* K(n) \sma_{\widehat{E(n)}} K(n)^{op}$ which gives rise to $\bar{q}_i$ in the spectral sequence converging to $\pi_* THH^{\widehat{E(n)}}(K(n))$ also lives in $\pi_* k(n) \sma_S k(n)^{op}$, under the name $\bar{\tau}_i$. The rest is a simple matter of comparing two ways to calculate $\pi_* THH^S(k(n))$ in low degrees, either by first running the B\"okstedt spectral sequence and then the Adams spectral sequence, or by running the K\"unneth spectral sequence.
\end{proof}

\begin{theorem} \label{Kntraceob}
The canonical map
\begin{equation}
K(n) \sar \Sigma^{2(p^n-1)/(p-1)-n-2p^i-1} \widehat{E(n)}/v_i
\end{equation}
extends to a $\mathcal{W}_{(n-i)(p-1)}$-trace and the obstruction to a $\mathcal{W}_{(n-i)(p-1)+1}$-trace is nontrivial, giving a height $(n-i)(p-1)$ extension
\begin{equation}
v_i \bar{q}_i^{p-1} \cdots \bar{q}_{n-1}^{p-1} = v_n
\end{equation}
in the spectral sequence converging to $\pi_* THH^{\widehat{E(n)}}(K(n))$.
\end{theorem}

Here $\bar{q}_i^{p-1}$ should be interpreted as $-\gamma_{p-1}(\bar{q}_i)$.

\begin{proof}
This follows from Theorem \ref{K(n)obstructions} and Proposition \ref{prop:sigmatau=qbar}, using that this is the first possible obstruction for degree reasons.
\end{proof}

\begin{corollary}
The canonical map
\begin{equation}
\hEn/v_i \lar K(n)
\end{equation}
extends to a $\mathcal{W}_{(n-i)(p-1)}$-cotrace and the obstruction to a $\mathcal{W}_{(n-i)(p-1)+1}$-cotrace is $(-1)^{n-i} v_n q_i^{p-1} \cdots q_{n-1}^{p-1}$.
\end{corollary}

This is enough to determine $THH(K(1))$, and almost enough to determine $THH(K(n))$ for $n>1$.

\begin{theorem} \label{thm:THHK(1)calc}
Topological Hochschild cohomology of $K(1)$ is given by
\begin{equation}
\pi_* THH_{\widehat{E(1)}}(K(1)) \cong \Z_p[v_1,v_1^{-1}][[q]]/(p+v_1q^{p-1})
\end{equation}
as a ring.

Similarly,
\begin{equation}
\pi_* THH^{\widehat{E(1)}}(K(1)) \cong \bigoplus_{i=0}^{p-2} \Sigma^{2i} \Z/p^\infty[v_1,v_1^{-1}]
\end{equation}
as an $\widehat{E(1)}_*$-module.
\end{theorem}

\begin{proof}
This follows from Theorem \ref{Kntraceob}, after we use the Weierstrass preparation theorem to conclude that $\Z_p[v_1,v_1^{-1}][[q]]/(p+v_1q^{p-1}+a_{2p-2} q^{2(p-1)}+\ldots)$ does not depend on the coefficients $a_{2p-2},\ldots$.
\end{proof}

Thus $\widehat{E(1)}_* \sar \pi_* THH_{\widehat{E(1)}}(K(1))$ is a tamely ramified (at $p$) extension of degree $p-1$, while the $p$-completion, or $K(1)$-localization, of $THH^{\widehat{E(1)}}(K(1))$ consists of $p-1$ copies of of $\widehat{E(1)}$.
It is a curious fact that both $THH_{\widehat{E(1)}}(K(1))$ and $KU^\wedge_p$ consist of $p-1$ copies of $\widehat{E(1)}$, though these copies are glued together in different ways.

\begin{remark}
While $THH_{\widehat{E(1)}}(K(1))$ is an $E_2$ ring spectrum by the Deligne conjecture, it is not an $E_\infty$ ring spectrum. One can see this by considering suitable power operations in $K(1)$-local $E_\infty$ ring spectra. Recall, e.g.~from \cite{Re04} that a $K(1)$-local $E_\infty$ ring spectrum $T$ (which has to satisfy a technical condition which we do not have to worry about here) has power operations $\psi$ and $\theta$ such that (in particular) $\psi$ is a ring homomorphism and
\begin{equation}
\psi(x)=x^p+p\theta(x)
\end{equation}
for $x \in T^0 X$. Now, if $T_*$ has an $i$'th root of some multiple of $p$, say, $\zeta^i=ap$ for a unit $a$ and $i>1$, then we get
\begin{equation}
ap=\psi(\zeta)^i=(\zeta^p+p\theta(\zeta))^i,
\end{equation}
and the right hand side is divisible by $p^2$ while the left hand side is not. In particular, we can apply this to $T=THH_{\widehat{E(1)}}(K(1))$ as above to show that this cannot be an $E_\infty$ ring spectrum.

We believe that a similar argument shows that the spectra $THH_{\hEn}(K(n))$ and $THH_{E_n}(K_n)$ can never be $E_\infty$, except in the cases when $THH_{E_n}(K_n) \simeq E_n$.
\end{remark}

In general we get the following.

\begin{theorem}
Topological Hochschild cohomology of $K(n)$ is given by
\begin{equation}
\pi_* THH_{\hEn}(K(n)) \cong \hEn_*[[q_0,\ldots,q_{n-1}]]/(p-f_0,\ldots,v_{n-1}-f_{n-1})
\end{equation}
as a ring, where $f_i \equiv (-1)^{n-i} v_n q_i^{p-1} \cdots q_{n-1}^{p-1}$ modulo higher degree.

Similarly,
\begin{equation}
\pi_* THH^{\hEn}(K(n)) \cong \Gamma_{\hEn_*}[\bar{q}_0,\ldots,\bar{q}_{n-1}]/(x_i \gamma_I(\bar{q})-\sum a_{iJ} \gamma_{I-J}(\bar{q})),
\end{equation}
where the notation is the same as in Corollary \ref{cor:allTHHforM}.
\end{theorem}

\begin{conjecture} \label{conj:THHofK(n)unique}
Topological Hochschild cohomology of $K(n)$ is independent of the $A_\infty$ structure and a finite, tamely ramified (at the primes $p, v_1,\ldots, \\ v_{n-1}$) extension of degree $(p-1) \cdots (p^n-1)$.
\end{conjecture}

One way to check this conjecture would be to calculate more of the coefficients of the power series $f_i$. We believe, using the philosophy that if something can happen it will, that they look like
\begin{eqnarray*}
f_{n-1} & = & -v_n q_{n-1}^{p-1}+\ldots \\
f_{n-2} & = & v_n q_{n-1}^{p-1} \pm v_n q_{n-2}^{p^2-1}+\ldots \\
\vdots & & \\
f_0 & = & (-1)^n q_0^{p-1} \cdots q_{n-1}^{p-1}
\pm v_n q_0^{p-1} \cdots q_{n-2}^{p^2-1} +\ldots \pm v_n q_0^{p^n-1}+\ldots
\end{eqnarray*}
in which case the conjecture would follow. The part about being independent of the $A_\infty$ structure would also follow (by using that $THH_{\hEn}(K(n)) \simeq THH_S(K(n))$) if we knew that all the $A_\infty$ structures on $K(n)$ become equivalent over $S$. We might come back to that elsewhere.

For the $2$-periodic Morava $K$-theories $K_n=E_n/(p,u_1,\ldots,u_{n-1})$ we have many more choices of $A_\infty$ structures, and hence more choices for $THH(K_n)$. In particular, we have the following.

\begin{theorem} \label{THHK_n1}
For any $p$ and $n$ there exists $A_\infty$ structures on $K_n$ such that
\begin{equation}
THH_{E_n}(K_n) \simeq E_n
\end{equation}
and
\begin{equation}
THH^{E_n}(K_n) \simeq E_n[I^{-1}]/IE_n.
\end{equation}
\end{theorem}

\begin{proof}
If $n>1$ or $p$ is odd, this follows from Theorem \ref{thm:Cinvertible}. The last case is Baker and Lazarev's calculation of $THH_{KU}(KU/2)$.
\end{proof}

In general we get many possible extensions, and while we can choose most of the coefficients of the power series $f_i$ freely, we run into the same problem as in the calculation of $THH(K(n))$. Again the case $n=1$ is easier than the general case.

\begin{theorem}
Given an $n$ with $1 \leq n<p-1$ and $a \in \{1,\ldots,n-1\}$, there is an $A_\infty$ structure on $K_1$ with
\begin{equation}
\pi_* THH_{E_1}(K_1) \cong (E_1)_*[[q]]/(p+a(uq)^n).
\end{equation}
For such an $A_\infty$ structure, $\pi_*THH^{E_1}(K_1)$ is a direct sum of $n$ copies of $\Z/p^\infty[u,u^{-1}]$. Otherwise,
\begin{equation}
\pi_* THH_{E_1}(K_1) \cong (E_1)_*[[q]]/(p+(uq)^{p-1}),
\end{equation}
and $\pi_* THH^{E_1}(K_1)$ is a direct sum of $p-1$ copies of $\Z/p^\infty[u,u^{-1}]$.
\end{theorem}

Thus $THH_{E_1}(K_1)$ is always a finite extension of $E_1$, of degree $d$ for some $1 \leq d \leq p-1$. As in Conjecture \ref{conj:THHofK(n)unique} we believe that $THH_{E_n}(K_n)$ is always a finite extension of $E_n$.

\section{$THH$ of Morava $K$-theory over $S$}
In this section, which is somewhat different from the previous sections, we prove that $THH(K(n))$ and $THH(K_n)$ do not depend on the ground ring. By this we mean that the canonical maps
\begin{equation}
THH^S(K_n) \lar THH^{E_n}(K_n)
\end{equation}
and
\begin{equation}
THH_{E_n}(K_n) \lar THH_S(K_n)
\end{equation}
are weak equivalences, and similarly for $K(n)$ using $\hEn$ instead of $E_n$. The earliest incarnation of this equivalence can be found in \cite{Ro88}, where Robinson observed that for $p$ odd the $t_i$'s in
\begin{equation}
\pi_*(K(n) \sma_S K(n)) \cong K(n)_*[\alpha_0,\ldots,\alpha_{n-1},t_1,t_2,\ldots]/(\alpha_i^2,v_n t_i^{p^n}-v_n^{p^i}t_i)
\end{equation}
do not contribute to the $Ext$ groups $Ext^{**}_{\pi_* K(n) \sma_S K(n)}(K(n)_*,K(n)_*)$. Something similar is true at $p=2$ if we use $K(n) \sma_S K(n)^{op}$. While $\alpha_i$ squares to $t_{i+1}$ instead of $0$ in this case, the $Ext$ calculation is still valid. This was used by Baker and Lazarev in \cite{BaLa} to see that $THH_S(KU/2) \simeq THH_{KU}(KU/2)$.

Much of the material in this section comes from \cite{Re97}, where Rezk does something similar to show that certain derived functors of derivations vanish. We have also used ideas from \cite{HoMi}.

We expect $THH$ to be invariant under change of ground ring from $S$ to $E_n$, or the other way around, because something similar holds algebraically.

\begin{lemma}
Let $R \lar R'$ be a Galois extension of rings and suppose $A$ is an $R'$ algebra. Then the canonical maps
\begin{equation}
HH_*^R(A) \lar HH_*^{R'}(A)
\end{equation}
and
\begin{equation}
HH^*_{R'}(A) \lar HH^*_R(A)
\end{equation}
are isomorphisms
\end{lemma}

\begin{proof}
Recall from \cite{GeWe} that Hochschild homology satisfies \'etale descent and Galois descent. \'Etale descent shows that $HH_*^{R'}(A) \cong HH_*^R(A)$ when $A=A' \otimes_R R'$, and then Galois descent shows that it holds for any $A$. The cohomology case is similar.
\end{proof}

Now, if $E_n$ is the Morava $E$-theory associated to the Honda formal group over $\F_{p^n}$, Rognes describes \cite[\S 5.4]{Ro_Galois} how the unit map $S \lar E_n$ is a $K(n)$-local (or $K_n$-local) pro-Galois extension with Galois group $\mathbb{G}_n$, the extended Morava stabilizer group. Similarly, $S \lar \hEn$ is a $K(n)$-local pro-Galois extension with the slightly smaller Galois group $\mathbb{G}_n/K$ for $K=\F_{p^n}^\times \times Gal(\F_{p^n}/\F_p)$ , so we expect the result, if not the proof, to carry over.

\subsection*{Perfect algebras}
Let $A$ and $B$ be commutative $\F_p$-algebras, and suppose $i:A \lar B$ is an algebra map. There is a Frobenius map $F$ sending $x$ to $x^p$ on each of these $\F_p$-algebras. Let $A^F$ denote $A$ regarded as an $A$-algebra using the Frobenius $F$. Now we can define a relative Frobenius $F_A: A^F
\otimes_A B \lar B$ as $F_A(a \otimes b)=i(a)b^p$ on decomposable tensor factors, or as the unique map $A^F \otimes_A B \sar B$ in the following diagram:
\begin{equation}
\xymatrix{
A \ar[r]^i \ar[d]_F & B \ar[d] \ar[rdd]^F & \\
A \ar[r] \ar[rrd]_i & A^F \otimes_A B \ar@{.>}[rd] & \\
 & & B }
\end{equation}

\begin{definition}
We say that $i:A \lar B$ is \emph{perfect} if $F_A:A^F \otimes_A B \lar B$ is an isomorphism.
\end{definition}

This definition specializes to the usual definition of a perfect $\F_p$-algebra when $A=\F_p$.

Now suppose that $i:A \lar B$ has an augmentation $\epsilon : B \lar A$. Let $I=ker(\epsilon)$ be the
augmentation ideal, so that $B \cong A \oplus I$ additively.

\begin{lemma}
For $i \geq 0$ and any $B$-module $M$ we have
\begin{equation}
Tor^B_i(I,M) \cong Tor^B_{i+1}(A,M)
\end{equation}
and
\begin{equation}
Ext^i_B(I,M) \cong Ext^{i+1}_B(A,M).
\end{equation}
\end{lemma}

\begin{proof}
This follows by choosing a resolution like
\begin{equation}
A \longleftarrow A \oplus I \longleftarrow P_0 \longleftarrow P_1 \longleftarrow \ldots,
\end{equation}
of $A$, where $P_0 \longleftarrow P_1 \longleftarrow \ldots$ is a projective resolution of $I$ as a $B$-module.
\end{proof}

Now, if $i:A \lar B$ is perfect, we have an isomorphism $F_A : A^F \otimes_A (A \oplus I) \lar A \oplus I$, and
this gives an isomorphism $F_A : A^F \otimes_A I \lar I$ of non-unital algebras.

Now suppose that $M$ is an $A$-module, and regard $M$ as a $B$-module via $\epsilon$. Then $I$ acts trivially on
$M$, and we use that to prove the following:

\begin{proposition}
Suppose that $i:A \lar B$ is perfect and let $M$ be any $A$-module viewed as a $B$-module via the augmentation $\epsilon : B \lar A$. Then we have
\begin{equation}
Tor^B_i(I,M)=0
\end{equation}
and
\begin{equation}
Ext^i_B(I,M)=0
\end{equation}
for all $i$.
\end{proposition}

\begin{proof}
We show that the maps
\begin{equation}
(A^F \otimes_A I) \otimes_B M \arrow{F_A \otimes 1} I \otimes_B M 
\end{equation}
and
\begin{equation}
Hom_B(I,M) \overset{F_A^*}\lar Hom_B(A^F \otimes_A I,M)
\end{equation}
are both isomorphisms and zero. They are isomorphism because $i: A \lar B$ is perfect. They are zero because, for example, given any map $f:I \lar M$ of $B$-modules, we find that $F_A^*f$ is given by $F_A^*f(a \otimes b)=f(ab^p)=b^{p-1}f(ab)=0$, so $F_A^*f$ is zero. The same argument applies to a projective resolution of $I$ to show that $Ext^i_B(I,M)=0$ for $i>0$. The argument for $Tor$ is similar.
\end{proof}

Combining the above two results we get the following:

\begin{theorem} \label{thm:perfect}
Suppose that $i:A \lar B$ is perfect, and let $M$ be any $A$-module regarded as a $B$-module via the
augmentation $\epsilon : B \lar A$. Then we get
\begin{equation}
Tor^B_i(A,M)=0
\end{equation}
and
\begin{equation}
Ext_B^i(A,M)=0
\end{equation}
for $i>0$, while $A \otimes_B M \cong M$ and $Hom_B(A,M) \cong M$.
\end{theorem}

\subsection*{Formal groups}
Most of what we need to know about formal groups can be found in \cite{Re97}. Recall that given two Morava $E$-theories $E$ and $F$ of the same height, the maximal ideals in $E_0 F$ coming from $m_E$ and $m_F$ coincide. Furthermore, $E_0F/m$ represents isomorphisms of formal group laws. Let $W(\Gamma_1,\Gamma_2)=k_1 \otimes_L W \otimes_L k_2$, where $L$ is the Lazard ring (isomorphic to $MU_*$, or $MUP_0$, where $MUP$ is the $2$-periodic complex cobordism spectrum) and $W=L[t_0^{\pm 1},t_1,\ldots]$.

\begin{proposition} (\cite[Remark 17.4]{Re97})
If $E$ and $F$ are the Morava $E$-theories associated to two formal groups $\Gamma_1$ and $\Gamma_2$ of height $n$, then
\begin{equation}
(\pi_0 E \sma_S F)/m \cong W(\Gamma_1,\Gamma_2).
\end{equation}
\end{proposition}

\begin{proposition} (\cite[Corollary 21.6]{Re97})
The ring $W(\Gamma_1, \Gamma_2)$ is a perfect $k_1$-algebra.
\end{proposition}

\begin{proposition} \label{structureofKsmaKop} $ $
Given any multiplication on $K$, we have
\begin{equation}
\pi_0(K \sma_S K^{op}) \cong \big(\pi_0(E \sma_S E)\big)/m \otimes \Lambda(\alpha_0,\ldots,\alpha_{n-1}).
\end{equation}
additively, and each $\alpha_i$ squares to something that acts trivially on $K_*$.
\end{proposition}

\begin{proof}
This is clear additively, and the claim about the multiplicative structure follows as in the proof of Proposition \ref{prop:AAopring}.
\end{proof}

Now we are in a position to prove the following theorem:

\begin{theorem} \label{THHInv}
Let $E$ be either $E_n$ or $\hEn$. If $E=E_n$ let $K=K_n$ and if $E=\hEn$ let $K=K(n)$. Then the canonical maps
\begin{equation}
THH^S(K) \lar THH^E(K)
\end{equation}
and
\begin{equation}
THH_E(K) \lar THH_S(K)
\end{equation}
are weak equivalences.
\end{theorem}

\begin{proof}
We have spectral sequences calculating $\pi_*$ of both sides, where the $E_2$-terms are $Tor_{\pi_*(K \sma_EK^{op})} (K_*,K_*)$ and $Tor_{\pi_*(K \sma_SK^{op})}(K_*,K_*)$ in the first case and the corresponding $Ext$ groups in the second case. For $(E,K)=(E_n,K_n)$, Theorem \ref{thm:perfect} and Proposition \ref{structureofKsmaKop} shows that the $E_2$-terms are isomorphic, and since the isomorphisms are induced by the obvious maps this proves the theorem.

The case $(E,K)=(\hEn,K(n))$ is similar, using $L[t_1,t_2,\ldots]$ instead of $W$.
\end{proof}

One interesting consequence of this theorem is the following:

\begin{corollary} \label{cor:samespaceofAinfty}
Let $E$ be either $E_n$ or $\hEn$. If $E=E_n$ let $K=K_n$ and if $E=\hEn$ let $K=K(n)$. Then the spaces of $A_\infty$ $E$-algebra structures on $K$ and $A_\infty$ $S$-algebra structures on $K$ are equivalent.
\end{corollary}

Note that this is before modding out by $Aut(A)$, which acts on the space of $A_\infty$ structures. Because $Aut_S(K)$ is larger than $Aut_E(K)$, there are fewer equivalence classes of $S$-algebra structures on $K$.

\bibliographystyle{plain}
\bibliography{b}
\vspace{12pt}
\noindent
Department of Mathematics, University of Chicago
\newline
5734 S University Ave
\newline
Chicago, IL 60637
\newline
Email: vigleik@math.uchicago.edu

\end{document}